%
%
\input amstex
\documentstyle{amsppt}


\def \ni {\noindent}

\def \m {\medskip}
\def \b {\bigskip}
\def \a {\alpha}
\def \A {\Bbb A}
\def \ad {{\text {ad}}}
\def \be {\beta}

\def \d {\delta}
\def\D {\Delta}

\def \g {{\frak g}}
\def \ga {\gamma}
\def \h {{\frak h}}

\def\isom{\cong}
\def\J {\Bbb J}
\def \l {\lambda}
\def\Lnum#1#2#3{\left[ \matrix #1\ ; \  #2 \\ \ \ #3\ \endmatrix \right]}
\def \o {\overline}
\def \ot {\otimes}
\def\qbin#1#2#3{\left[ \matrix #1 \\ #2 \endmatrix \right]_{#3}}

\def \C {\Bbb C}
\def \Z {\Bbb Z}




\topmatter
\title {Verma-type Modules for} \\
{Quantum Affine Lie Algebras}\endtitle
\leftheadtext{FUTORNY, GRISHKOV, MELVILLE}
\author
Vyacheslav M. Futorny\\
Alexander. N. Grishkov \\
Duncan J. Melville
\endauthor
\rightheadtext{QUANTUM VERMA-TYPE MODULES}

\address Instituto de Matematica, Universidade do Sao Paulo,
Sao Paulo, Brasil.
\endaddress
\email futorny\@ime.usp.br
\endemail
\address Instituto de Matematica, Universidade do Sao Paulo,
Sao Paulo, Brasil.
\endaddress
\email grishkov\@ime.usp.br
\endemail
\address Department of Mathematics,
St. Lawrence University,
Canton, New York  13617   USA.
\endaddress
\email dmel\@music.stlawu.edu
\endemail


\abstract
Let $\g$ be an untwisted affine Kac-Moody algebra and $M_J(\l)$
a Verma-type module for $\g$ with $J$-highest weight
$\l \in P$.  We construct quantum Verma-type modules
$M_J^q(\l)$ over the quantum group $U_q(\g)$, investigate
their properties and show that $M_J^q(\l)$ is a true quantum
deformation of $M_J(\l)$ in the sense that the weight structure
is preserved under the deformation.  We also analyze the submodule
structure of quantum Verma-type modules.
\endabstract

\thanks The first author is a Regular Associate of the ICTP.\endgraf
The third author was supported in part by a Faculty Research Grant from
St. Lawrence University.
\endthanks

\date March 14, 2000\enddate




\endtopmatter

\document


\NoBlackBoxes

\subhead Introduction \endsubhead
\b
The representation theory of Kac-Moody algebras is much richer
than that of finite-dimensional simple algebras.  In particular,
Kac-Moody algebras have modules containing both finite and
infinite-dimensional weight spaces, something that cannot happen
in the finite-dimensional setting \cite{Le}. These representations
of Kac-Moody algebras arise from taking non-standard partitions of
the root system, partitions which are not Weyl-equivalent to the
standard partition into positive and negative roots.  For affine
algebras, there is always a finite number of Weyl-equivalency
classes of these nonstandard partitions.  Corresponding to each
partition is  a Borel subalgebra and one can form representations
induced from one-dimensional modules for these nonstandard Borel
subalgebras.  These induced modules are called modules of
Verma-type, in analogy with standard Verma modules induced from a
standard Borel subalgebra.

Verma-type modules were first studied and classified by Jakobsen
and Kac \cite{JK1, JK2}, and by Futorny \cite{Fu1, Fu2}. Further
work elucidating their structure, including the construction of
the appropriate categorical setting, determination of
irreducibility criteria, submodule structure, BGG duality and BGG
resolutions can be found in \cite{Co1, Co2, CFM, Fu4, Fu5, and
references therein}.

The theory of Verma-type modules is best developed in the case
where the central element of $\g$ acts with non-zero charge.  In
this case we say the level is non-zero.  Verma-type modules of
zero level are still not completely classified.  There are many
technical difficulties, for example, a Verma-type module of level
zero may have subquotients that are not quotients of Verma-type
modules \cite{Fu3, Proposition 5(ii)}.  As the affine theory is
not complete, we restrict our attention on structural questions to
representations of  non-zero level, except in the case of
imaginary Verma modules, where the level zero theory is more
complete.
 However, the techniques we use can be generalized (with
more technicalities) to the case of level zero where $\l$ is in
``general position'' (see \cite{Fu4} for details of the affine
theory).

Since their introduction by Drinfeld \cite{Dr} and Jimbo \cite{Ji}
in 1985, there has been a tremendous interest in studying quantized
enveloping algebras for Kac-Moody algebras and their representations.
Quantum groups have turned out to be extremely important objects with
rich and diverse connections to an ever-increasing number of areas
 of mathematics and physics.  A vigorous body of research is
developing on determining their structure, their representations
and their applications.

In many cases, the representation theory of quantum groups
parallels that of the associated underlying classical algebras,
although often with some subtle differences.  The closest parallel
comes when the classical and quantum representations have the same
weight structure. In this paper, we construct quantum Verma-type
modules over the quantum group $U_q(\g)$ associated to an
untwisted affine Kac-Moody algebra $\g$, prove they have the same
weight structure as the corresponding Verma-type modules, and
investigate some of their properties.

One of the problems to be faced when studying nonstandard
representations of a Kac-Moody algebra $\g$, is that the absence
of a general PBW theorem for quantum groups means that we cannot
lift the triangular decomposition of $\g$ up to a triangular
decomposition of $U_q(\g)$.  In \cite{CFKM}, the authors studied
the simplest non-standard case, that of Verma-type modules for the
algebra $U_q(A_1^{(1)})$.  In the case of $U_q(A_1^{(1)})$, there
is (up to Weyl-equivalence) only one non-standard partition and
the associated representations are called imaginary Verma modules.
In \cite{CFKM}, the authors had to construct an appropriate PBW
basis to deal with that particular case (see their Proposition
2.2).  The techniques used there do not easily generalize to the
case of all affine algebras.

In this paper we rely heavily on the work by Beck \cite{Be1, Be2}
and Beck and Kac \cite{BK} on PBW bases for quantum groups of
affine algebras, and we exploit a particularly convenient
description of the nonstandard partition of the root system used
to construct the quantum  Verma-type modules.  This approach leads
to one of our key ideas, that it is possible to determine a basis
for the quantum  Verma-type modules in a unified manner without
giving a PBW basis for the algebra.

Having constructed the quantum representations, we  check that
they are true quantum deformations of the equivalent classical
modules.  That is, the quantum and classical modules have the same
weight structure.  To do this, we follow the $\A$-form technique
originally introduced by Lusztig \cite{Lu}, and subsequently
refined and developed by Kang and co-authors \cite{Ka, CFKM,
BKMe}.  For an overview of this procedure and a summary of known
quantum deformation results, see \cite{M}.  Our main result of
this section, generalizing that of \cite{CFKM}, is that any
quantum Verma-type module with integral $J$-highest weight $\l$ is
a quantum deformation of the equivalent Verma-type module over
$U(\g)$ for $\g$ an untwisted affine Kac-Moody algebra.

Using the quantum deformation theorem, we study some of the
structural properties of quantum Verma-type  modules.  In
particular, we prove a general irreducibility criterion for
quantum Verma-type modules and probe the structure of quantum
imaginary Verma modules at level zero, when they are reducible.
The results obtained are similar to those given for non-quantum
imaginary Verma modules in \cite{Fu3}, showing that these quantum
modules are closely related to their classical cousins. Some of
the results on quantum imaginary Verma modules obtained in this
paper were announced in \cite{FGM}.

The structure of the paper is as follows.  First, we recall
background information and establish notation in Section 1.
In Section 2 we review the construction of Verma-type modules for
affine algebras.  In Section 3 we construct quantum
Verma-type modules and provide them with a basis.  Section
4 constructs $\A$-forms, and Section 5 gives the classical
limits and quantum deformation theorem.  Section 6 discusses the
structural results for Verma-type modules in general and Section 7
considers quantum imaginary Verma modules at level zero.
For additional basic background material and notation on Kac-Moody
algebras, see the book by Kac \cite{K}; for background
information on quantum groups, see the excellent texts by
Chari and Pressley \cite{CP} and Jantzen \cite{Ja}.

\b
\subhead 1. Preliminaries \endsubhead
\b

\ni {\bf 1.1.} Let $N$ be a positive integer.  Fix index sets
$\dot I = \{ 1, \dots, N\}$ and $I = \{ 0, \dots, N\}$. Let
$\dot \g$ be a finite-dimensional simple complex Lie algebra with
Cartan subalgebra $\dot \h$, root system
$\dot \D \subset \dot \h^*$,
and set of simple roots $\dot \Pi = \{\a_1,\dots,\a_N\}$. Denote
by $\dot \D_+ $ and $\dot \D_-$, the positive and negative  roots
of $\g$. Let $\dot Q = \oplus_{i=1}^N \Z \a_i$ be the root lattice
of $\dot \g$, and let $\dot A = (a_{ij})_{1\le i,j \le N}$ be the
Cartan matrix for $\dot \g$.  Define a basis $h_1,\dots, h_N$ of
$\dot \h$ by $\a_i(h_j) = a_{ij}$.  Let $\dot P=\{\l \in \h^*\ |\
\l(h_i)\in \Z, i=1,\dots ,N \}$ be the weight lattice of $\dot
\g$. Let $(.|.)$ denote both the symmetric invariant bilinear form
on $\dot \g$ and the induced form on $\dot \g^*$, normalized so
that $(\a|\a)=2$ for any short root $\a$.  For $i = 1,\dots, N$,
let $d_i = (\a_i|\a_i)/2$.  Then each $d_i$ is a positive integer,
the $d_i$ are relatively prime, and the diagonal matrix
$\dot D = \text{diag}(d_1,\dots, d_N)$ is such that $\dot D \dot A$
is symmetric.

\m
\ni {\bf 1.2.} Let $\g$ denote
the untwisted affine Kac-Moody algebra associated to $\dot \g$.
Then $\g$ has the loop space realization
$$
\g = \dot \g \ot \C[t,t^{-1}] \oplus \C c \oplus \C d,
$$
where $c$ is central in $\g$; $d$ is the degree derivation,
so that $[d,x \ot t^n] = n x \ot t^n$ for any $x \in \dot \g$ and
$n \in \Z$, and $[x \ot t^n, y \ot t^m] = [x,y] \ot t^{n+m} +
\d_{n+m,0}n(x|y)c$ for all $x,y \in \dot \g$, $n,m \in \Z$.  We
set $\h = \dot \h \oplus \C c \oplus \C d$.

The algebra $\g$ has a Cartan matrix
$A = (a_{ij})_{0 \le i,j \le N}$ which is an extension of $\dot A$.
There exists an integer $d_0$ and a diagonal matrix
$D = \text{diag}(d_0, \dots, d_N)$ such that $DA$ is symmetric.
An alternative Chevalley-Serre presentation of $\g$ is given
by defining it as the Lie algebra with generators $e_i, f_i, h_i$
 ($i \in I$) and $d$ subject to the relations
$$
\aligned
[h_i,h_j] &= 0, \qquad [d,h_i]=0,\\
[h_i, e_j] &= a_{ij}e_j, \qquad [d,e_j] = \delta_{0,j}e_j,\\
[h_i, f_j] &= -a_{ij}f_j, \qquad [d, f_j] = -\delta_{0,j}f_j,\\
[e_i, f_j] &= \delta_{ij}h_i, \\
(\ad e_i)^{1-a_{ij}}(e_j)&= 0, \qquad
(\ad f_i)^{1-a_{ij}}(f_j) = 0, \quad i \neq j.
\endaligned
$$

\m
\ni {\bf 1.3.} We can define the root system of $\g$ in the following
way.
Extend the root lattice $\dot Q$ of $\dot \g$ to a lattice
$Q = \dot Q \oplus \Z\d$, and extend the form $(.|.)$ to $Q$
by setting $(q|\d)=0$ for all $q \in \dot Q$ and $(\d|\d)=0$.
The root system $\D$ of $\g$ is given by
$$
\D = \{ \a + n \d\ |\ \a \in \dot \D, n \in \Z\}
\cup
\{k\d\ |\ k\in \Z, k \neq 0\}.
$$
The roots of the form $\a+n\d$, $\a \in \dot \D, n \in \Z$ are
called real roots, and those of the form $k\d$,
$k \in \Z, k \neq 0$ are called imaginary roots.  We let
$\D^{re}$ and $\D^{im}$ denote the sets of real and
imaginary roots, respectively.  The set of positive real roots
of  $\g$ is $\D_+^{re} = \dot \D_+
\cup \{\a + n\d \ |\ \a \in \dot \D, n>0\}$ and the set of positive
imaginary roots is $\D_+^{im} = \{k\d\ |\ k>0\}$.
The set of positive roots of $\g$ is
$\D_+ = \D_+^{re} \cup \D_+^{im}$.  Let $Q_+$ be the monoid
generated by $\D_+$.
Similarly, on the negative
side, we have $\D_- = \D_-^{re}\cup \D_-^{im}$, where
$\D_-^{re} = \dot \D_- \cup \{ \a + n\d\ |\ \a \in \dot \D, n<0\}$
and $\D_-^{im} = \{k\d\ |\ k<0\}$ and let $Q_-$ be the monoid
generated by $\D_-$.
Further, if $\theta$ denotes the highest positive root of
$\dot \g$ and $\a_0 := \d - \theta$, then
$\Pi = \{\a_0,\a_1,\dots, \a_N\}$ is a set of simple
roots for $\g$.
We extend the weight lattice $\dot P$ of $\dot \g$ to the
weight lattice $P$ of $\g$ defined as
$P = \{ \l \in \h^*\ |\ \l(h_i) \in \Z, i \in I, \l(d) \in \Z\}$.
Let $W$ denote the Weyl group of
$\g$ generated by the simple reflections $r_0, r_1, \dots, r_N$
and $B$ denote the associated braid group with generators
$T_0, T_1, \dots, T_N$.

\m
\ni {\bf 1.4.}
Beck \cite{Be1, Be2} has introduced a total ordering of the root
system leading to  PBW bases for $\g$ and its quantum analog,
$U_q(\g)$.  We state the construction here, partially following the
more abstract notation developed by Damiani \cite{Da} and
Gavarini \cite{Ga}.  For a related PBW construction, see
\cite{KT}.

For any affine algebra $\g$, there exists a map $\pi :\Z \mapsto I$
such that, if we define
$$
\beta_k =
\left \{
\aligned
&r_{\pi(0)}r_{\pi(-1)}\cdots r_{\pi(k+1)}(\a_{\pi(k)})
\qquad \text{ for all } k \le 0 \\
&r_{\pi(1)}r_{\pi(2)} \cdots r_{\pi(k-1)}(\a_{\pi(k)})
\qquad  \text{ for all } k\ge1,
\endaligned
\right.
$$
then the map $\pi':\Z \mapsto \D_+^{re}$ given by
$\pi'(k)=\beta_k$
is a bijection.  Further, we can choose $\pi$ so that
$ \{ \beta_k\ |\ k \le 0\} =
\{ \a + n\d \ |\ \a \in \dot \D_+, n \ge 0\}$ and
$\{ \beta_k\ |\ k \ge 1\} =
\{ -\a + n\d\ |\ \a \in \dot \D^+, n>0\}$.

It will be convenient for us to invert Beck's original ordering
of the positive roots (cf. \cite{BK, 1.4.1} for the original
order and \cite{Ga, \S 2.1} for this ordering).  Thus, we set
$$
\beta_0 > \beta_{-1} > \beta_{-2}> \dots >\d >2\d>  \dots
>\beta_2 >\beta_1.
$$
Clearly, if we say $-\a < -\beta$ if and only if $\beta > \a$ for
all
positive roots $\a, \beta$, we obtain a corresponding ordering
on $\D_-$.

The following elementary observation on the ordering will play
a crucial role later.  Write $A<B$ for two sets $A$ and $B$ if
$x < y$ for all $x \in A$ and $y \in B$.  Then Beck's total
ordering of the positive roots can be divided into three sets:
$$
\{ \a + n\d\ |\ \a \in \dot \D_+, n \ge 0\} >
\{k\d \ |\ k>0\} > \{-\a+k\d\ |\ \a \in \dot \D_+, k>0\}.
$$
Similarly, for the negative roots, we have,
$$
\{ -\a - n\d\ |\ \a \in \dot \D_+, n \ge 0\} <
\{-k\d \ |\ k>0\} < \{\a-k\d\ |\ \a \in \dot \D_+, k>0\}.
$$

Note that the map $\pi$, and so the total ordering, is not unique.
We assume a suitable $\pi$ chosen and fixed now throughout the
paper.  Beck's original approach and proof is constructive, but the
existential approach avoids some technicalities we do not need
below.

\m
\ni {\bf 1.5.}
The quantum group, or quantized universal enveloping algebra,
of $\g$ is the associative algebra $U_q(\g)$ with 1 over $\C(q)$
with generators $E_i, F_i, K_i^{\pm 1}$ ($i\in I$) and
$D^{\pm 1}$  subject to the defining relations
$$
\aligned
&K_iK_i^{-1}= K_i^{-1}K_i = DD^{-1}=D^{-1}D=1,\\
&K_iK_j= K_jK_i, \qquad K_iD = DK_i, \\
&K_iE_j = q_i^{a_{ij}} E_jK_i, \qquad DE_j=q_0^{\d_{j,0}}E_jD, \\
&K_iF_j = q_i^{-a_{ij}} F_iK_i, \qquad DF_j= q_0^{-\d_{j,0}}F_jD,
\\
& E_iF_j - F_jE_i =
\d_{ij}\frac{K_i-K_i^{-1}}{q_i-q_i^{-1}},\\
& \sum_{k=1}^{1-a_{ij}}(-1)^k
\qbin{1-a_{ij}}{k}{q_i} E_i^{1-a_{ij}-k}E_jE_i^k = 0,
\qquad i \neq j \\
& \sum_{k=1}^{1-a_{ij}}(-1)^k
\qbin{1-a_{ij}}{k}{q_i} F_i^{1-a_{ij}-k}F_jF_i^k = 0,
\qquad i \neq j,
\endaligned
$$
where $q_i = q^{d_i}$ (we can choose $d_i$ so that $d_0 =1$ and
$q_0 = q$), and
$$
\qbin{m}{n}{q} = \frac{[m]_q!}{[m-n]_q![n]_q!},\qquad
[m]_q!= \prod_{j=1}^m [j]_q, \qquad
[j]_q = \frac{q^j-q^{-j}}{q-q^{-1}}
$$
for all $i \in I$, $m,n \in \Z$, $m\ge n >0$.
For any $\mu \in Q$, we have $\mu = \sum_{i\in I} c_i\a_i$,
for some integers $c_i$.  Denote $K_{\mu} = \prod_{i\in I}
K_i^{c_i}$.  Then $K_{\l}K_{\mu}=K_{\l+\mu}$ for all
$\l,\mu \in Q$.  In particular, we have
$K_{\pm \a_i} = K_i^{\pm 1}$.
Let $U_q^+(\g)$ (resp. $U_q^-(\g)$) be the subalgebra of
$U_q(\g)$ generated by $E_i$ (resp. $F_i$), $i \in I$, and
let $U_q^0(\g)$ denote the subalgebra generated by
$K_i^{\pm}$ ($i \in I$) and $D^{\pm}$.

The action of the braid group generators $T_i$ on the generators
of the quantum group $U_q(\g)$ is given by the following.
$$
\aligned
T_i(E_i) &= -F_i K_i, \qquad T_i(F_i) = -K_i^{-1}E_i, \\
T_i(E_j) &= \sum_{r=0}^{-a_{ij}} (-1)^{r-a_{ij}}
\frac{1}{[-a_{ij}-r]_{q_i}!}\frac{1}{[-r]_{q_i}!} q_i^{-r}
E_i^{-a_{ij}-r}E_jE_i^{r}, \qquad \text{if } i \neq j,\\
T_i(F_j) &=  \sum_{r=0}^{-a_{ij}} (-1)^{r-a_{ij}}
\frac{1}{[-r]_{q_i}!}\frac{1}{[-a_{ij}-r]_{q_i}!} q_i^r
F_i^{r}F_jF_i^{-a_{ij}-r}, \qquad \text{if } i \neq j,\\
T_i(K_j) &= K_jK_i^{-a_{ij}}, \qquad
T_i(K_j^{-1}) = K_j^{-1}K_i^{a_{ij}}, \\
T_i(D) &= DK_i^{-\delta_{i,0}}, \qquad
T_i(D^{-1}) = D^{-1}K_i^{\delta_{i,0}}.
\endaligned
$$

\m
\ni{\bf 1.6.}
For each $\beta_k \in \D_+^{re}$, define the root vector
$E_{\beta_k}$ in $U_q(\g)$ by
$$
E_{\beta_k} =
\left \{
\aligned
&E_{\pi(0)} \qquad k=0\\
&T^{-1}_{\pi(0)}T^{-1}_{\pi(-1)}\cdots
T^{-1}_{\pi(k+1)}(E_{\pi(k)})
\qquad \text{ for all } k < 0 \\
&E_{\pi(1)} \qquad k = 1 \\
&T_{\pi(1)}T_{\pi(2)} \cdots T_{\pi(k-1)}(E_{\pi(k)})
\qquad  \text{ for all } k>1.
\endaligned
\right.
$$
Each real root space is 1-dimensional, but each imaginary root
space is $N$-dimensional.  Hence, for each positive imaginary
root $k\d$ ($k>0$) we define $N$ imaginary root vectors,
$E_{k\d}^{(i)}$ ($i\in \dot I$) by
$$
\exp\left(
(q^i-q^{-i})\sum_{k=1}^{\infty}E_{k\d}^{(i)}z^k\right) =
1 + (q^i-q^{-i})\sum_{k=1}^{\infty} K_i^{-1}[E_i,
E_{-\a_i+k\d}]z^k.
$$
Then for each $k$, the $E_{k\d}^{(i)}$ span the $k\d$-root space
and commute with each other.  Further, the $E_{\beta_k}$ ($k\in
\Z$)
and $E_{k\d}^{(i)}$ ($k>0$) form a basis for $U_q^+(\g)$.

Let $\omega$ denote the standard $\C$-linear antiautomorphism of
$U_q(\g)$, and set $E_{-\a} = \omega(E_{\a})$ for all $\a \in
\D_+$.
Then $U_q$ has a basis of elements of the form $E_-HE_+$,
where $E_{\pm}$ are ordered monomials in the $E_{\a}$,
$\a \in \D_{\pm}$, and $H$ is a monomial in $K_i^{\pm 1}$
and $D^{\pm 1}$ (which all commute).

Furthermore, this basis is, in Beck's terminology, convex, meaning
that, if $\a, \beta \in \D_+$ and $E_{\beta} > E_{\a}$, that is,
$\beta > \a$, then
$$
E_{\beta}E_{\a} - q^{(\a|\beta)}E_{\a}E_{\beta} =
\sum_{\a<\ga_1<\dots <\ga_r <\be}c_{\ga}
E_{\ga_1}^{a_1}\cdots E_{\ga_r}^{a_r}
$$
for some integers $a_1,\dots, a_r$ and scalars
$c_{\ga}\in \C[q, q^{-1}]$,
$\ga = (\ga_1,\dots,\ga_r)$ \cite{BK, Proposition 1.7c}, and
similarly for the negative roots.

\b
\subhead 2. Verma-type modules for affine algebras \endsubhead
\b

We recall here the construction and properties of Verma-type
modules for affine algebras.  The theory of Verma-type modules has
been developed in a number of papers, for example
\cite{Co1, Co2, CFM, Fu2, and Fu4}. For convenience, we
use \cite{Fu5} as our
standard reference. It contains detailed proofs and references to
the original publication of the results summarized below.

The root system $\D$ of $\g$ has a natural partition into positive
and negative roots, $\D^+$ and $\D^-$.  An arbitrary partition
$\D = S \cup -S$ is called {\it closed} if whenever
$\a$ and $\beta$ are in $S$ and $\a + \be$ is a root, then
$\a + \be \in S$.  If $S$ is a closed partition, then the space
$\g_S = \oplus_{\a \in S}\g_{\a}$ is a subalgebra of $\g$, and
$\g$ has a triangular decomposition
$\g = \g_{-S} \oplus \h \oplus \g_S$.

Let $W$ denote the
Weyl group of $\g$.  Then it is well-known that, for any
finite-dimensional complex simple Lie algebra, all closed
partitions are $W\times \{\pm 1\}$-equivalent to the standard
partition into positive and negative roots.  For Kac-Moody
algebras this is no longer true.  In particular, for affine
algebras, there are always a finite number (greater than 1)
of inequivalent Weyl-orbits of partitions.  These non-standard
partitions were first studied and classified
by Jakobsen and Kac \cite{JK1, JK2} and Futorny
\cite{Fu1, Fu2}.

Next we introduce some notation we will need dealing with
the root system of $\g$ and construct a collection
of sets that parametrize closed partitions.

Let $J \subseteq \dot I = \{1, \dots, N\}$.  Let
$\Pi^J = \{ \a_j \in \Pi\ |\ j \in J\}$.  Set
$Q^J = \oplus_{j \in J} \Z \a_j \oplus \Z\delta$, and
$Q_{\pm}^J = Q^J \cap Q_{\pm}$.  For any
subset $\epsilon$ of $\dot \Pi$, let
$Q^{\epsilon}_{\pm}$ denote a semigroup of $\h^*$ generated
by $\pm \epsilon$.

Let $\dot \D^J$ be the finite root system generated by the
simple roots in $\Pi^J$.  Then $\dot \D^J = \emptyset$ if
and only if $J = \emptyset$.  Set
$$
\D^J = \{ \a + n\d \in \D\ |\ \a \in \dot \D^J, n \in \Z\}
\cup \{n\d \ |\ n \in \Z \setminus \{0\} \},
$$
and let $\D_{\pm}^J = \D^J \cap \D_{\pm}$ and
$\D_J^{\infty} = \D \setminus \D^J$.
Now we let
$$
{_+\D_J^{\infty}} = \{ \a + n\d \in \D\ |\
\a \in \dot \D_+ \setminus \dot \D^J, n \in \Z \}.
$$
Finally, let $S_J = \D^J_+ \cup {_+\D_J^{\infty}}$.

Note that, for any $J \subseteq \dot I$, the positive
imaginary roots $k\d$ ($k>0$) are in $\D^J_+$ and
a simple root $\a_i \in \D^J_+$ if and only if $i \in J$,
while $\a_i \in {_+\D_J^{\infty}}$ if and only if $i \in \dot I
\setminus J$.  Hence, regardless of choice
of $J$, we will always have $\dot \D_+ \cup \{ k\d\ |\ k>0\}
\subseteq S_J$.  Further, we note that $S_J$ is a closed
partition of $\D$ for any $J \subseteq \dot I$.  The two
extreme cases are when $J = \dot I$, in which case
$S_J = \D_+$, and $J = \emptyset$, in which case
$$
S_{\emptyset} = \{ \a + n\d\ |\ \a \in \dot \D_+, n \in \Z\}
\cup \{k \d\ |\ k>0\}.
$$
The sets $S_J$ parametrize the closed partitions of $\D$.  That
is, a closed partition $S$ is $W \times \{\pm 1\}$-equivalent
to a unique $S_J$ for some $J \subseteq \dot I$
\cite{Fu5, Theorem 2.4}.

Verma-type modules are representations induced from
1-dimensional representations for nonstandard Borel subalgebras.
Equivalence classes of Borel subalgebras correspond to
equivalence classes of closed partitions.  In order to
classify Verma-type modules, we begin by classifying the
Borel subalgebras for an affine algebra $\g$.  We will
always assume the fixed Cartan subalgebra $\h$.

Let $\sigma$ denote the Chevalley involution of $\g$
defined by $\sigma(e_i) = -f_i, \sigma(f_i) =-e_i,
\sigma(h) = -h$ for $i \in I, h \in \h$.  Let
$B \subset \g$ be a Lie subalgebra such that
$\h \subset B$ and $B + \sigma(B)
= \g$.  Then we call $B$ a Borel subalgebra
if $B \cap \sigma(B) = \emptyset$, and
a parabolic subalgebra otherwise.  The algebra
involution $\sigma$ induces a linear involutive
antiautomorphism $\sigma$ on the root system $\D$,
whereby $\sigma(\a)= -\a$ for all $\a \in \D$.  Hence,
there is a correspondence between Borel subalgebras
of $\g$ and closed partitions of $\D$.

Let $U(\g_S)$ (resp. $U(\g_{-S})$) denote the universal enveloping
algebra of $\g_S$ (resp. $\g_{-S}$).  Then, by the PBW theorem,
the triangular decomposition of $\g$ determines a triangular
decomposition of $U(\g)$ as
$U(\g) = U(\g_{-S}) \ot U(\h) \ot U(\g_S)$.

Let $\l \in \h^*$.  Then $\l$ extends to a map on $(U(\h))^*$, also
denoted by $\l$.
A $U(\g)$-module $V$ is called a weight module
if $V= \oplus_{\mu \in \h^*} V_{\mu}$, where
$V_{\mu} = \{v \in V\ |\ h \cdot v = \mu(h) v \text{ for all } h
\in U(\h)\}$.  The non-zero
subspaces $V_{\mu}$ are called weight spaces.  Any submodule
of a weight module is a weight module.

Let $B$ be a Borel subalgebra and $\l:B\to \C$ a 1-dimensional
representation.  A weight $U(\g)$-module $V$ is said to be of
{\it highest weight} $\l$ with respect to $B$ if there is some
nonzero vector $v \in V$ such that $V = U(\g)\cdot v$ and
$x \cdot v = \l(x)v$ for all $x \in B$.  We define the
induced module
$$
M_B(\l) = U(\g) \otimes_{U(B)} \C,
$$
called the module of {\it Verma-type} associated with $B$ and $\l$.
If $B$ is the standard Borel subalgebra of $\g$, then $M_B(\l)$ is
just the usual Verma module with highest weight $\l$.  If $B$ is
not $W\times\{\pm 1\}$-conjugate to the standard Borel, then
the Verma-type module $M_B(\l)$ contains both finite and
infinite-dimensional weight spaces.

Since there is a 1-to-1 correspondence between Borel
subalgebras and closed partitions of the root system,
and the sets $S_J$ parametrize the Weyl-equivalence
classes of partitions, we can construct a set of
representatives of the Borel conjugacy classes and so a
canonical collection of Verma-type modules.  For
$J \subseteq \dot I$, let $B^J = \sum_{\beta \in S_J} \g_{\beta}
\oplus \h$.  Then $B^J$ is a Borel subalgebra of $\g$.
Denote $M_{B^J}(\l)$ by $M_J(\l)$.  Then the $M_J(\l)$,
$J \subseteq \dot I$, parametrize all Verma-type modules.
We call $M_J(\l)$ a $J$-highest-weight module with $J$-highest
weight $\l$.
Henceforth, we will consider only the Verma-type modules
$M_J(\l)$.

The {\it level} of $M_J(\l)$ is $\l(c)$.  At certain points,
we need to make the assumption that
$\l(c) \neq 0$. The case of level zero is more
complicated and not completely understood.  Although some
of our techniques do generalize in certain case, we cannot
give complete results.

Recalling the two special cases of $J$, we have that if $J = \dot I$,
then $S_J = \D_+$ and $M_{\dot I}(\l)$ is the ordinary Verma module
of highest weight $\l$.  At the other extreme, the modules
$M_{\emptyset}(\l)$ are called {\it imaginary Verma modules}
\cite{Fu3}.

We denote by $\frak{n}_{\pm J}$ the subalgebras $\g_{\pm S_J}$, so we
have a triangular decomposition $\g = \frak{n}_{-J} \oplus \h
\oplus \frak{n}_{J}$.

In the proposition below we collect some basic statements about
the structure of the Verma-type modules $M_J(\l)$.  Proofs
can be found in \cite{Fu5, Propositions 3.4 and 5.2}.

\proclaim{Proposition 2.1}
Let $\l \in \h^*$ and $J \subseteq \dot I$, and let $M_J(\l)$
be the Verma-type module
of $J$-highest weight $\l$.  Then $M_J(\l)$ has the following properties.
\roster
\item"(i)" The module $M_J(\l)$ is a free
$U(\frak{n}_{-J})$-module of rank 1
generated by the $S$-highest weight vector $1 \ot 1$ of weight $\l$.
\item"(ii)" $M_J(\l)$ has a unique maximal submodule and hence
a unique irreducible quotient, which we denote $L_J(\l)$.
\item"(iii)" Let $V$ be a $U(\g)$-module generated by some
$J$-highest weight vector $v$ of weight $\l$.  Then there exists a
unique surjective homomorphism $\phi:M_J(\l) \mapsto V$ such
that $\phi(1\ot1) = v$.
\item"(iv)" $\dim M_J(\l)_{\mu} \neq 0$ if and only if $\l-\mu$ is in
the monoid generated by $S_J$, \newline
$\dim M_J(\l)_{\l} = 1$, \newline
and $0< \dim M_J(\l)_{\mu} < \infty$ if and only if $\l -\mu \in Q^J_+$.
\endroster
\endproclaim

Let $\g_J := \sum_{\beta \in \D^J} \g_{\beta} \oplus \h$.
We call a subset $C \subseteq J$ connected if the
Coxeter-Dynkin diagram associated to the simple roots
$\a_i$, $i\in C$, is connected.  The set $J$ can then be
partitioned into a collection $\Cal C$ of subsets corresponding
to connected components of the Coxeter-Dynkin diagram associated
to $J$.  For $C \in \Cal C$, let $\frak{A}_C$ be the affine
subalgebra of $\g$ generated by $e_i, f_i$ ($i \in C$), together
with the central element $c$ and degree derivation $d$ from $\g$.
Let $\g^f = \sum_{C \in \Cal C} \frak{A}_C$ and
$\tilde \g^f = \g^f + \h$.

Let $G:= \oplus_{k \in \Z \setminus \{0\}} \g_{k\d} \oplus \C c$.  The
algebra $G$ is called a Heisenberg subalgebra of $\g$.
$G$ has a triangular decomposition $G = G_- \oplus \C c \oplus G_+$,
where $G_{\pm} = \oplus_{k>0} \g_{\pm k \d}$.
Set $\o G = \{ g \in G \ |\ [g,[\g^f,\g^f]]=0\}$,
and set $\o G_{\pm} = \o G \cap G_{\pm}$.
By Theorem 3.3 of
\cite{Fu5}, we have the following proposition.

\proclaim{Proposition 2.2}
$\g_J = \o G_-\oplus \tilde \g^f \oplus \o G_+$.
\endproclaim

By Proposition 2.1 (iv), we know that
$0 < \dim M_J(\l)_{\mu} < \infty$ if and only if
$ \l-\mu \in Q_+^J$.  Let us say
$\mu \le_J \l$ if $\l -\mu \in Q_+^J$.
Let $M^f_J(\l) := \oplus_{\mu \le_J \l} M_J(\l)_{\mu}$.
Then $M^f_J(\l)$ is the sum of all the finite-dimensional
weight spaces of $M_J(\l)$.  We also set
$L_J^f(\l)= \oplus_{\mu \le \l}L_J(\l)_{\mu}$.

\proclaim{Proposition 2.3 \rm(\cite{Fu5, Proposition 5.3})}
$M_J^f(\l) \isom \tilde M^f(\l) \otimes \o M(\l)$, where
$\tilde M^f(\l)$ is a Verma module for $\tilde \g^f$ and
$\o M(\l)$ is a Verma module for $\o G$.
\endproclaim

Next we recall some statements about the submodule structure of
$M_J(\l)$.  Up to this point, we have not needed the
restriction that $\l(c) \neq 0$, but we do need it for
parts of the next proposition.
Let $\frak{u}_{\pm J} = \oplus_{\beta \in S_J \setminus \D^J}
\g_{\pm \beta}$.

\proclaim{Proposition 2.4 \rm(\cite{Fu5, Lemma 5.4, Theorem 5.14})}
Assume
 $\l \in \h^*$.
\roster
\item"(i)" If $0 \neq v \in M_J(\l)$, then
$U(\g)v \cap M_J^f(\l) \neq 0$.
\item"(ii)" If $N$ is a $\g$-submodule of $M_J(\l)$ and $\l(c) \neq 0$, then
$$ N \isom U(\frak{u}_{-J}) \ot_{\C} ( N \cap M_J^f(\l))
$$
as vector spaces.
\item"(iii)" If $\l(c) \neq 0$, then $L_J(\l) \isom U(\frak{u}_{-J})\ot_{\C} L_J^f(\l)$
as vector spaces.
\endroster
\endproclaim

We have $\l \in \h^*$.  For $C \in \Cal C$, let $\l_C = \l|_{(\h
\cap \frak{A}_C)^*}$, and let $M^C(\l)$ denote the Verma module
over $\frak{A}_C$ with highest weight $\l_C$.  Since $\tilde
M^f(\l)$ is isomorphic to a tensor product of $M^C(\l)$, we obtain
from Proposition 2.4 the following irreducibility criterion.

\proclaim{Corollary 2.5}  The Verma-type module $M_J(\l)$
is irreducible if and only if $\l(c) \neq 0$ and
$M^C(\l)$ is an irreducible $\frak{A}_C$-module for
every $C \in \Cal C$.
\endproclaim

In particular, for the special case of imaginary Verma modules,
we have the following simple criterion.

\proclaim{Corollary 2.6 \rm(\cite{Fu5, Proposition 5.8})}  The
imaginary Verma module $M_{\emptyset}(\l)$ is irreducible
if and only if $\l(c) \neq 0$.
\endproclaim

\b
\subhead 3. Verma-type modules for quantum affine algebras
\endsubhead
\b

Next we must construct the quantum modules that we will show are
deformations of the classical Verma-type modules.  Let $J
\subseteq \dot I$. Let $\frak{n}^q_{\pm J}$ be the
subalgebra of $U_q(\g)$ generated by $\{ E_{\beta}\ |\ \beta \in
\pm S_J\}$, and let $B_q$ denote the subalgebra of $U_q(\g)$
generated by $\{ E_{\beta}\ |\ \beta \in S_J\}\cup U_q^0$.

A $U_q(\g)$-module
$V^q$ is called a quantum weight module if
$V^q = \oplus_{\mu \in P}V_{\mu}^q$, where
$$
V_{\mu}^q = \{ v \in V\ |\ K_i^{\pm 1}\cdot v = q_i^{\pm \mu(h_i)}v,
D^{\pm 1}\cdot v = q_0^{\pm \mu(d)}v \}.
$$
Any submodule of a quantum weight module is a weight module.
A $U_q(\g)$-module $V^q$ is called a $J$-highest weight
module with highest weight $\l\in P$ if there is a non-zero vector
$v \in V^q$ such that:
\roster
\item"(i)" $u^+\cdot v  =0$ for all
$u^+ \in \frak{n}^q_{ J} \setminus \C(q)^*$;
\item"(ii)" for each $i \in I$,
$K_i^{\pm 1}\cdot v = q_i^{\pm\l(h_i)}v$,
$D^{\pm 1}\cdot v = q_0^{\pm \l(d)}v$;
\item"(iii)" $V^q = U_q(\g) \cdot v$.
\endroster
Note that, in the absence of a general quantum PBW theorem for
non-standard partitions, we cannot immediately claim that  a
$J$-highest weight module $V^q$ is generated by
$\frak{n}^q_{- J}$.
This is in contrast to the classical case, and the reason behind
Theorem 3.5 below.

We define a $U_q(\g)$-module as follows. Let $\C(q) \cdot v$
be a 1-dimensional vector space.
Let $\l \in P$, and set  $E_{\beta}\cdot v=0$ for all
$\beta \in S_J$,  $K_i^{\pm 1}\cdot v = q_i^{\pm\l(h_i)}v$
($i \in I$) and
$D^{\pm 1}\cdot v = q_0^{\pm \l(d)}v$.
Now define $M_J^q(\l)=U_q(\g)\ot_{B_q} \C(q) v$.  Then
$M_J^q(\l)$ is a $J$-highest weight $U_q$-module called
the {\it quantum Verma-type module} with $J$-highest weight $\l$.

The following contains the basic properties of Verma-type modules.

\proclaim{Proposition 3.1}
Let $\l \in P$ and $J \subseteq \dot I$, and let $M^q_J(\l)$
be the quantum Verma-type module
of $J$-highest weight $\l$.  Then $M^q_J(\l)$ has the following
properties.
\roster
\item"(i)" The module $M^q_J(\l)$ is a free
$\frak{n}^q_{- J}$-module of rank 1
generated by the $J$-highest weight vector $v_{\l}$ of weight $\l$.
\item"(ii)" $M^q_J(\l)$ has a unique maximal submodule and hence
a unique irreducible quotient, which we denote $L^q_J(\l)$.
\item"(iii)" Let $V$ be a $U_q(\g)$-module generated by some
$J$-highest weight vector $v$ of weight $\l$.  Then there exists a
unique surjective homomorphism $\phi:M^q_J(\l) \mapsto V$ such
that $\phi(v_{\l}) = v$.
\endroster
\endproclaim

We also want to show that the module $M_J^q(\l)$
 is spanned by the
``right'' set of vectors.  In order to do this, we must
appeal to Beck's PBW basis of $U_q(\g)$
and to the very useful grading by degree introduced by Beck and
Kac \cite{BK, \S 1.8}, which we reproduce here as our notation
differs slightly from theirs.

Beck \cite{Be2} has shown that $U_q(\g)$ has a basis comprising
elements of the form $N_{(a_{\beta})}KM_{(a'_{\beta})}$,
where the $M_{(a_{\beta})}$ are ordered monomials in
$E_{\beta}^{a_{\beta}}$, $\beta \in \D_+$, $a_{\beta} \in \Z_+$,
$N_{(a_{\beta})} = \omega(M_{(a_{\beta})})$, and $K$ is an
ordered monomial in $K_i^{\pm}$ and $D^{\pm}$.  The
notation $(a_{\beta})$ indicates the sequence of powers $a_{\beta}$
as $\beta$ runs over $\D_+$.  Of course, almost all terms of the
sequence are zero.

In \cite{BK, \S 1.8}, Beck and Kac define the {\it total height} of
such a basis element by
$$
d_0(N_{(a_{\beta})}KM_{(a'_{\beta})}) = \sum_{\beta \in \D_+}
(a_{\beta}+a'_{\beta}) \text{ht}\beta,
$$
where $\text{ht} \beta$ is the usual height of a root.  Next, they
set the {\it total degree} of a basis element to be
$$
d(N_{(a_{\beta})}KM_{(a'_{\beta})}) =
(d_0(N_{(a_{\beta})}KM_{(a'_{\beta})}),(a_{\beta}),(a'_{\beta}))
\in \Z_+^{2\D_+ +1}.
$$
Considering $\Z_+^{2\D_+ +1}$ as  a totally ordered semigroup
with the usual lexicographical ordering, Beck and Kac introduce a
filtration of $U_q(\g)$ by defining $U_s$, for any
$s \in \Z_+^{2\D_+ +1}$,
to be the span of the basis monomials
$N_{(a_{\beta})}KM_{(a'_{\beta})}$ with
degree $d(N_{(a_{\beta})}KM_{(a'_{\beta})}) \le s$.  Finally,
they obtain the following proposition.

\proclaim{Proposition 3.2 \rm(\cite{BK, Proposition 1.8})}
The associated graded algebra $\text{Gr} U_q(\g)$ of the
$\Z_+^{2 \D_+ +1}$-filtered algebra $U_q(\g)$ is the algebra
over $\C(q)$ generated by $E_{\a}$, $\a \in \D$, counting
multiplicities, $K_i^{\pm 1}$ ($i \in I$), and $D^{\pm}$
 subject to the relations
$$
\aligned
&K_iK_i^{-1}= K_i^{-1}K_i = DD^{-1}=D^{-1}D=1,\\
&K_iK_j= K_jK_i, \qquad K_iD = DK_i, \\
&K_iE_{\a} = q^{(\a_i |\a)} E_{\a}K_i,
\qquad DE_{\a}=q^{n}E_{\a}D, \text{ for } \a = \gamma + n\delta,
\gamma \in \dot \Delta. \\
&E_{\a}E_{-\beta} = E_{-\beta}E_{\a} \text{ if } \a, \beta \in \D_+,\\
&E_{\a}E_{\beta} = q^{(\a| \beta)}E_{\beta}E_{\a}, \qquad
E_{-\a}E_{-\beta} = q^{(\a|\beta)}E_{-\beta}E_{-\a},
\text{ if } \a, \beta \in \D_+ \text{ and }\beta < \a.
\endaligned
$$
\endproclaim

\m
Next we need the following technical result, which we state in an
abstract manner.  Let $\Cal I$ be a totally ordered set without
infinite decreasing chains (i.e., $\Cal I$ has a minimal element with
respect to the ordering).  Let $A$ be an associative algebra over
a field $K$ with generators $\{ a_i \ |\ i \in \Cal I\}$.  Let
$O=\{(n_i|i\in \Cal I, n_i\in \Z_+, n_i=0
\text{ for all but finite number of indices })\}$.  The set $O$ has
a total lexicographical order
such that $(n_i\ |\ i\in \Cal I)>(m_i\ |\ i\in \Cal I)$ if and
only if there exists
some $j\in \Cal I$ such that $n_j>m_j$ and
$n_k=m_k$ for all $k>j$.

Suppose that the algebra $A$ has a basis
$B=\{v=a_{i_1}^{k_1}\dots a_{i_s}^{k_s}\ |\ i_1>\dots>i_s\}$.
Then the orderings of $\Cal I$ and $O$ impose an ordering on $B$.
For any word $v=a_{j_1}...a_{j_k}\in A$, denote by $\o{v}$ the
unique element in $B$ such that
$v=\lambda \o{v}+ (\text{ terms lower in the ordering })$, for some
$\lambda\in K$.  If
$\o{v}=a_{i_1}^{k_{i_1}}\dots a_{i_s}^{k_{i_s}}$ as an element of B,
set $|v|=(k_{i_1},\dots ,k_{i_s})\in O$.

\proclaim{Proposition 3.3} Suppose that for all $i$ and $j$ in
$\Cal I$ we have:
(*) $a_ia_j=\xi_{ij}a_ja_i +\sum_{v\in B}\xi_{v}v$, where
$\xi_{ij}\neq 0$ and
$|a_ia_j|>|v|$ for all $v$ such that $\xi_{v}\neq 0$. Then
$|vw|=|v|+|w|$ for all $v, w\in B$.
\endproclaim
\demo{Proof} We proceed by induction on $|w|$ and, for fixed
$|w|$, by induction on $|v|$. By our assumption the set $\Cal I$
has a minimal element $i_0$.
If $|w|=\epsilon_{i_0}=(n_i|n_i=0, i\neq i_0, n_{i_0}=1)$
then $w$ is the minimal element in $O$.  Hence, for any
$v \in B$, we have  $\o{vw}=vw$ and so $|vw|=|\o{vw}|=
|v|+|w|$.

Let $w=a_i$ and $v=v_1a_j$.  If $j\geq i$,  then
$\o{vw}=vw$ and $|vw|=|v|+|w|$. If  $j<i$, then by (*),
$vw=v_1a_ja_i=\xi_{ij}v_1 a_ia_j +\sum_{u\in B}\xi_{u}v_1u$.
We will prove that $|v_1u|<|v_1a_ia_j|$. By induction on $|v|$ we
have $|v_1a_i|=|v_1|+|a_i|$, while by induction on $|w|$ we get  $|v_1a_ia_j|=|v
_1a_i|+|a_j|= |v_1|+|a_i|+|a_j|$.
Let $\xi_{u} \neq 0$ and $u=a_ku_1$. If $k<i$ then by induction on
$|w|$ we have $|v_1u|=|v_1|+|u|$, since $|u|<|a_i|$.
 If $k=i$ then $u_1=a_su_2$ and $s<j<i$ since
$|u|=|a_ia_su_2|<|a_ia_j|$. Then
$|v_1u|=|v_1a_ia_su_2|=|v_1a_i|+|a_su_2|$ and, as $|a_i|>|a_su_2|$,
we also have
$|v_1a_i|+|a_su_2| =|v_1|+|a_i|+ |a_su_2|<|v_1|+|a_i|+|a_j|$,
since $s<j$.
Hence we have $|v_1u|<|v_1a_ia_j|$ if $\xi_{u}\neq 0$ and
$|vw|=|v_1a_ia_j|=|v_1|+|a_i|+|a_j|=|v|+|a_i|$.

Now assume that $w=w_1a_i$ with $|w|>|w_1|>|a_i|$. By induction,
$|vw|=|vw_1a_i|=|vw_1|+|a_i|=|v|+|w_1|+|a_i|=|v|+|w|$. Hence, the
proposition is proved.
\hskip 1cm \qed
\enddemo

\m
Let $u$ be an arbitrary element of $U_q(\g)$.  We may write $u$
uniquely as a sum of basis monomials and define the total degree
of $u$ to be the largest total degree of these basis elements.  We
make the following observation.

\proclaim{Proposition 3.4} Let $u \in U_q(\g)$ be an arbitrary element
and $u', u''$ two basis monomials with $d(u') < d(u'')$.  Then
$d(uu')< d(uu'')$ and $d(u'u)<d(u''u)$.
\endproclaim
\demo{Proof}
The result follows from Propositions 3.2 and 3.3.
\hskip 1cm \qed
\enddemo

Next, we give the main result of this section.  First, however, we
introduce some notation to clarify the argument.  We will need these
definitions again in Section 6.

Consider the following subsets of $\D$:
$$
\aligned
A_1 &= \{ \a + n\d\ |\ \a \in \dot \D_+, n \ge 0\}, \\
A_2 &= \{ k\d\ |\ k>0\}, \\
A_3 &= \{-\a+k\d\ |\ \a \in \dot \D_+, k>0\}, \\
B_1 &= \{-\a -n\d\ |\ \a \in \dot \D_+, n \ge 0\},\\
B_2 &= \{ -k\d\ |\ k>0\}, \\
B_3 &= \{\a-k\d\ |\ \a \in \dot \D_+, k>0\}.
\endaligned
$$
Then $\D_+ = A_1\cup A_2\cup A_3$ and
$\D_- = B_1 \cup B_2 \cup B_3$.  Note that, in our ordering
of the root system, we have
$$
B_1<B_2<B_3<A_3<A_2<A_1.
$$

Now we must split these subsets into those associated with
finite and infinite-dimensional subspaces.  For $i = 1, 2, 3$,
let $A_i^{fin} = A_i \cap \D^J$,
$A_i^{\infty} = A_i \cap \D_J^{\infty}$,
$B_i^{fin} = B_i \cap \D^J$ and
$B_i^{\infty} = B_i \cap \D_J^{\infty}$.
Now, let
 $X_i$ denote an ordered monomial in
elements $E_{\beta}$, $\beta \in A_i$, and
$Y_i$ denote an ordered monomial of elements $E_{\beta}$,
$\beta \in B_i$.  Further, let
$X_i^{fin}$ (resp. $X_i^{\infty}$) denote an ordered monomial
in $E_{\beta}$, $\beta \in A_i^{fin}$ (resp. $\beta \in A_i^{\infty}$),
and let
$Y_i^{fin}$ (resp. $Y_i^{\infty}$) denote an ordered monomial in
$E_{\beta}$, $\beta \in B_i^{fin}$ (resp. $\beta \in B_i^{\infty}$).
We note that the sets $A_2^{\infty}$ and $B_2^{\infty}$ are
actually empty.

\proclaim{Theorem 3.5} As a vector space over $\C(q)$,
$M_J^q(\l)$ is isomorphic to
the space spanned by the ordered monomials
$$
E_{-\beta-n\d}\dots E_{-\beta+k\d}\dots
E_{-\a-n\d}\dots E_{-k\d}\dots E_{\a-k\d},
$$
for
$\a \in \dot \D_+\cap \D^J$, $\beta \in \dot \D_+ \cap \D_J^{\infty}$,
$n\ge 0$, $k>0$.
\endproclaim
\demo{Proof}
Utilizing Beck's PBW basis, we know any
element $u \in U_q(\g)$ can be written in the form
$$
u = \sum Y_1Y_2Y_3ZX_3X_2X_1,
$$
where $Z \in U_q^0(\g)$.

Let $v_{\l}$ be the canonical generator of $M_J^q(\l)$.  Suppose
$w \in M_J^q(\l)$.  Then, since $M_J^q(\l) = U_q(\g) \cdot v_{\l}$, we have
$w = u\cdot v_{\l}$ for some $u \in U_q(\g)$.  In view of the discussion
above, we may write $w = \sum Y_1Y_2Y_3ZX_3X_2X_1\cdot v_{\l}$,
for suitable monomials $X_i, Y_i, Z$.

As observed in Section 2, all roots of the form $\a + n\delta$ and
$k \delta$ ($\a \in \dot\D_+, n \ge 0, k >0$) are in $S_J$
for all subsets $J \subseteq \dot I$.  Hence,
monomials of the form $X_1$ and
$X_2$ all act as 0 on $v_{\l}$.  Further, $Z$ commutes with $X_3$ up
to a scalar in $\C(q)$.  Hence, we can write
$w = \sum Y_1Y_2Y_3X_3\cdot v_{\l}$.  The theorem asserts
that $M_J^q(\l)$ is spanned by monomials of the form
$Y_1^{\infty}X_3^{\infty}Y_1^{fin}Y_2Y_3^{fin}\cdot v_{\l}$.  We proceed
to rearrange the monomials in steps, arguing by induction on the
total degree.  The base of the induction is trivial as,
if $d_0=1$, there is only one simple root involved and nothing
to do.  Now suppose the theorem is true for all monomials
up to total degree $d$.  We consider monomials of next
highest degree.

First, we show that we can reduce to considering expressions
of the form $Y_1Y_2Y_3X_3^{\infty}\cdot v_{\l}$.  To do this, we
have to show that any monomial of the form $X_3^{fin}$ acts as
0 on $v_{\l}$, and that we can move all terms in $X_3^{fin}$ to the
right of any terms from $X_3^{\infty}$.

Since $A_3 \subset \D^+$, if a root $-\a +k \d$ is in
$A_3^{fin}= A_3 \cap \D^J$, then $-a +k\d \in \D_+^J$ and so
$-a+k\d\in S_J$.  Thus, $E_{-a+k\d}\cdot v_{\l}=0$, and we have
$X_3^{fin}\cdot v_{\l} = 0$.

Suppose a factor of the form $E_{-a+k\d}E_{-\beta+l\d}$ with
$-\a+k\d \in A_3^{fin}$ and $-\beta+l\d \in A_3^{\infty}$
occurs in a monomial $X_3$.  Since $E_{-\a+k\d}$ and $E_{-\beta+l\d}$
are both positive root vectors, we can use the convexity properties
of Beck's basis and the grading of the algebra to see that
$$
E_{-\a+k\d}E_{-\beta+l\d} = q^{\pm(\a | \beta)}E_{-\beta+l\d}E_{-a+k\d}
+fE_{-(\a+\beta)+(k+l)\d} + \sum f_{\gamma}E_{\gamma},
$$
where $\deg E_{\gamma} < \deg E_{-(\a+\beta)+(k+l)\d}$ and $f$ and
$f_{\gamma}$ are scalars in $\C[q,q^{-1}]$. The sign of
$q^{\pm (\a | \beta)}$ depends on whether $-a + k\d > -\beta+l\d$ or not.

The monomial $X_3'$ obtained from $X_3$ by replacing the factor
$E_{-a+k\d}E_{-\beta+l\d}$ by $E_{-\beta+l\d}E_{-\a+k\d}$ has
the $X_3^{fin}$ term moved to the right.  If $-a+k\d \in A_3^{fin}$
and $-\beta+l\d \in A_3^{\infty}$, then
$-(\a+\beta) + (k+l)\d \in A_3^{\infty}$.  Hence, the monomial $X_3''$
with $E_{-\a+k\d}E_{-\beta+l\d}$ replaced by
$E_{-(\a+\beta)+(k+l)\d}$ has the same total degree but one fewer term
in $X_3^{fin}$. We may repeat the process as necessary.
Finally, the monomials $X_3^{\gamma}$ with
$E_{-a+k\d}E_{-\beta+l\d}$ replaced by $E_{\gamma}$ are of
lower total degree by Proposition 3.4,
and can be rearranged by the inductive hypothesis.
Therefore, we need only consider the action of monomials of the
form $Y_1Y_2Y_3X_3^{\infty}$.

Secondly, we must determine how to commute monomials of the forms
$Y_3$ and $X_3^{\infty}$.   Let
$$
\aligned
X_3^{\infty} &= E_{-\a_1+k_1\d}\dots E_{-\a_r+k_r\d}, \text { and }\\
Y_3 &= E_{\beta_1-m_1\d}\dots E_{\beta_s-m_s\d},
\endaligned
$$
for suitable roots $-\a_i + k_i\d \in A_3^{\infty}$ ($i=1,\dots, r$)
and $\beta_j-m_j\d \in B_3$ ($j=1, \dots, s$).
Then we can write
$$
\aligned
Y_3X_3^{\infty} &=
E_{\beta_1-m_1\d}\dots E_{\beta_s-m_s\d}
E_{-\a_1+k_1\d}\dots E_{-\a_r+k_r\d}\\
&= E_{\beta_1-m_1\d}\dots E_{\beta_{s-1}-m_{s-1}\d}
E_{-\a_1+k_1\d}E_{\beta_s-m_s\d}E_{-\a_2+k_2\d}
\dots E_{-\a_r+k_r\d}\\
&\qquad + (\text{terms of lower total degree})
\endaligned
$$
by using the grading of Proposition 3.2.  Repeating this process
we get that
$$
Y_3X_3^{\infty}\cdot v_{\l} = X_3^{\infty}Y_3\cdot v_{\l} +
(\text{terms of lower degree} ) \cdot v_{\l}.
$$
By induction on the total degree, we may order the
terms of lower degree (as they act on $v_{\l}$), and we have
reduced the monomials we are concerned about to those
of the form $Y_1Y_2X_3^{\infty}Y_3$.

For the next step in the reduction process, we have to
show that any monomial of the form $Y_3^{\infty}$ acts
as 0 on $v_{\l}$ and that in $Y_3$ we can move all terms in
$Y_3^{\infty}$ to the right of any terms in $Y_3^{fin}$.

Suppose $Y_3^{\infty}$ contains a root vector $E_{\beta-k\d}$.
Then $\beta-k\d \in A_3 \cap \D_J^{\infty}$.  Since
$\beta-k\d \in \D_J^{\infty}$, we have $\beta-k\d \notin \D^J$ and
so $\beta \notin \dot \D^J$.  But since $\beta-k\d \in A_3$,
$\beta \in \dot \D_+$ and we must have
$\beta \in \dot \D_+\setminus \dot \D^J$.  Thus,
$\beta - k\d \in {_+\D_J^{\infty}} \subset S_J$ and
$E_{\beta-k\d}\cdot v_{\l} =0$.  Hence,
$Y_3^{\infty} \cdot v_{\l}=0$.

Now suppose a factor of the form
$E_{\beta-k\d}E_{\a-l\d}$, with $\beta-k\d \in B_3^{\infty}$ and
$\a-l\d \in B_3^{fin}$ occurs in a monomial $Y_3$.  Then both
$E_{\beta-k\d}$ and $E_{\a-l\d}$ are negative roots and we can
proceed with a similar argument to that used for showing we can
move factors in $X_3^{fin}$ to the right of those in $X_3^{\infty}$,
subject to our inductive hypothesis.  Note that in this case, if
$\beta-k\d\in B_3^{\infty}$ and $\a-l\d\in B_3^{fin}$, then
$(\a+\beta)-(k+l)\d \in B_3^{\infty}$ and the $Y_3^{\infty}$
factor absorbs the $Y_3^{fin}$ factor.

We are now reduced to rearranging monomials of the form
$Y_1Y_2X_3^{\infty}Y_3^{fin}$.  The fourth step is to
show that we can move monomials of the form $Y_2$ to the right
of those of the form $X_3^{\infty}$.

Consider a monomial $m$ containing a factor
$E_{-k\d}E_{-\a+l\d}$, $k>0$, $-\a+l\d \in A_3^{\infty}$.  Then
$E_{-k\d}$ is a negative root vector and $E_{-\a+l\d}$ is a
positive root vector and in the graded algebra, positive and
negative root vectors commute.  Thus, if $m'$ is the monomial
obtained from $m$ by replacing the factor $E_{-k\d}E_{-\a+l\d}$
by $E_{-\a+l\d}E_{-k\d}$, then we have
$$
m = m' + (\text{ terms of lower total degree }).
$$
The new terms of lower total degree may be disarranged, but, by
the inductive hypothesis, can be given the desired ordering.

The next step is to show that we can write monomials
$Y_1$ in the form $Y_1^{\infty}Y_1^{fin}$.  All root
vectors in $Y_1$ are negative and the argument is
similar to that used in the case of $X_3$ and $Y_1$.

The final step is to show that, up to terms of
lower degree, we can move monomials of the form
$Y_1^{fin}$ to the right of those of the form
$X_3^{\infty}$.  The argument here is similar to that
of the fourth step of the argument when we moved
monomials of the form $Y_2$ to the right of $X_3^{\infty}$.

It is important to note in this argument that, first,
we are rearranging monomials as they act on $v_{\l}$.  We do
not claim that these monomials are equivalent in the
algebra $U_q(\g)$.  That is, we are not claiming a PBW result.
Secondly, the ``commutations'' we have considered take place
with respect to the grading.  These monomials do not commute
in the algebra, or even directly with respect to their
action on $v_{\l}$.  At each stage, we may acquire more monomials,
but these additional terms can be reordered by the
inductive hypothesis.
\hskip 1 cm \qed
\enddemo

\b
\subhead 4. $\A$-forms of Verma-type modules \endsubhead
\b

In the previous section, we constructed quantum Verma-type
modules.  Now we show that these quantum Verma-type modules
are quantum deformations of Verma-type modules defined over
the affine algebra.  To do this, we need to show that the weight-space
structure of a given module $M_J^q(\l)$ is the same as that of its
classical
counterpart $M_J(\l)$ for any $\l \in P$ and $J \subseteq \dot I$.
The first step is to construct
an intermediate module, called an $\A$-form.

Following \cite{Lu}, for each $i \in I$, $s \in \Z$ and $n \in \Z_+$,
we define the {\it Lusztig elements} in $U_q(\g)$:
$$
\aligned
\Lnum{K_i}{s}{n} &=
\prod_{r=1}^n
\frac{K_iq_i^{s-r+1} - K_i^{-1}q_i^{-(s-r+1)}}{q_i^r-q_i^{-r}},\\
\Lnum{D}{s}{n} &= \prod_{r=1}^n
\frac{Dq_0^{s-r+1} - D^{-1}q_0^{-(s-r+1)}}{q_0^r-q_0^{-r}}.
\endaligned
$$

Let $\A= \C[q,q^{-1}, \frac{1}{[n]_{q_i}}, i\in I, n>0]$.  Define
the $\A$-form, $U_{\A}(\g)$, of $U_q(\g)$ to be the $\A$-subalgebra
of $U_q(\g)$ with 1 generated by the elements
$E_i, F_i, K_i^{\pm 1}, \Lnum{K_i}{0}{1}$, $i\in I$,
$D^{\pm 1}, \Lnum{D}{0}{1}$.
Let $U_{\A}^+$ (resp. $U_{\A}^-$) denote the subalgebra of
$U_{\A}$ generated by the $E_i$, (resp. $F_i$), $i \in I$, and
let $U_{\A}^0$ denote the subalgebra of $U_{\A}$ generated
by the elements $K_i^{\pm}, \Lnum{K_i}{0}{1}$,  $i \in I$,
$ D^{\pm 1}, \Lnum{D}{0}{1}$.

For any $i \in I$, $s\in \Z$ and $n \in \Z_+$, we have the following
identity
$$
\Lnum{K_i}{s}{n} = \prod_{i=1}^n \frac{1}{[r]_{q_i}}
\left( \Lnum{K_i}{0}{1} + [s-r+1]_{q_i}K_i^{-1}\right) \qquad
(\text{ cf. \cite{BKMe, eq. 3.8}}).
$$
Hence, all $\Lnum{K_i}{s}{n}$ are in $U_{\A}$.  Similarly,
all $\Lnum{D}{s}{n}$ are also in $U_{\A}$.

\smallskip
\proclaim{Proposition 4.1}
The following commutation relations hold between the generators
of $U_{\A}$.  For $i,j \in I$, $s \in \Z$, $n \in \Z_+$,
$$
\aligned
E_i \Lnum{K_j}{s}{n} & = \Lnum{K_j}{s-a_{ij}}{n}E_i, \\
\Lnum{K_j}{s}{n}F_i &= F_i \Lnum{K_j}{s-a_{ij}}{n}, \\
E_i \Lnum{D}{s}{n} &= \Lnum{D}{s-\delta_{i,0}}{n}E_i, \\
\Lnum{D}{s}{n}F_i &= F_i \Lnum{D}{s-\delta_{i,0}}{n},\\
E_iF_j &= F_j E_i, \qquad \text{ for } i \neq j, \\
E_iF_i^n &= F_i^nE_i + F_i^{n-1}\sum_{r=0}^{n-1}\Lnum{K_i}{-2r}{1}.
\endaligned
$$
\endproclaim
\demo{Proof} The first five equalities follow from the defining
relations of $U_q(\g)$ and the definition of the Lusztig elements,
while the last equality is proved by induction.
\hskip 1cm \qed
\enddemo

An immediate consequence of Proposition 4.1 is that $U_{\A}$
inherits the standard triangular decomposition of $U_q(\g)$.  In particular,
any element $u $ of $U_{\A}$ can be written as a sum of
monomials of the form $u^-u^0u^+$ where $u^{\pm} \in U_{\A}^{\pm}$
and $u^0 \in U_{\A}^0$.  In fact, we can say rather more.  For
each positive real root $\beta$, the root vector $E_{\beta}$ in
Beck's basis is defined via the action of the braid group on
the generators $E_i$.  But the coefficients of this action are all in
the ring $\A$.  Consequently, the real root vectors are in $U_{\A}$.
Next, consider the definition of the positive imaginary root vectors
$E_{k\delta}^{(i)}$, $i\in I$, $k>0$.  These are given
in terms of an exponential generating function containing commutators
of the form $[E_i, E_{-\a_i+k\delta}]$, and these will also
be in $U_{\A}$ since all the $E_i$ and $E_{-\a_i + k\delta}$ are.
The $\C(q)$ coefficients of the generating function
are all in $\A$, and so the imaginary root vectors are all in $U_{\A}$.
Thus, $U_{\A}$ inherits from $U_q(\g)$ a basis of monomials
of the form $N_{(a_{\beta})}KM_{(a'_{\beta})}$, where
$M_{(a'_{\beta})}$ and $N_{(a_{\beta})}$ are as before, and
$K$ is now an (ordered) monomial in the generators $K^{\pm}_i$,
$\Lnum{K_i}{0}{1}$, $i \in I$, $D^{\pm}$ and
$\Lnum{D}{0}{1}$ of $U_{\A}^{0}$.

\m
Let $\l \in P$, $J\subseteq \dot I$ and let $M_J^q(\l)$
be the Verma-type module
over $U_q(\g)$ with $J$-highest weight $\l$ and highest weight
vector $v_{\l}$.  The $\A$-form of  $M_J^q(\l)$, $M_J^{\A}(\l)$,
 is defined to be
the $U_{\A}$ submodule of $M_J^q(\l)$ generated by $v_{\l}$.
That is, we set
$$
M_J^{\A}(\l)= U_{\A} \cdot v_{\l}.
$$

\proclaim{Proposition 4.2} As a vector space over $\A$,
$M_J^{\A}(\l)$ is isomorphic
to the space spanned by the ordered monomials
$$
E_{-\beta-n\d}\dots E_{-\beta+k\d}\dots
E_{-\a-n\d}\dots E_{-k\d}\dots E_{\a-k\d},
$$
for
$\a \in \dot \D_+\cap \D^J$, $\beta \in \dot \D_+ \cap \D_J^{\infty}$,
$n\ge 0$, $k>0$.
\endproclaim
\demo{Proof}As in the proof of Theorem 3.5, we note that any
element $u$  in $U_{\A}$ can be written as a sum of monomials
of the form $Y_1Y_2Y_3ZX_3X_2X_1$, where the $X_i$ and
$Y_i$ are as in the theorem and $Z$ is now in $U_{\A}^0$.
Let $w \in M_J^{\A}(\l)$.  Then $w = u \cdot v_{\l}$ for
some $u \in U_{\A}$.  Write
$w = \sum Y_1Y_2Y_3ZX_3X_2X_1\cdot v_{\l}$.  As before,
we have $X_1 \cdot v_{\l} = 0$, and $X_2 \cdot v_{\l}=0$.  Also,
$Z$ commutes with $X_3$, up to a scalar in $\A$, by Proposition 4.1.

Now we must check the action of $Z$ on $v_{\l}$.  First, we
have $K_i^{\pm} \cdot v_{\l} = q_i^{\pm \l(h_i)} v_{\l} \in \A v_{\l}$,
and $D^{\pm} \cdot v_{\l} = q_0^{\pm \l(d)} v_{\l} \in \A v_{\l}$.
It remains only to check the action of the Lusztig elements.  For
$i \in I$, $s \in \Z$ and $n \in \Z_+$, we have
$$
\Lnum{K_i}{s}{n}\cdot v_{\l} =
\qbin{\l(h_i) + s}{n}{q_i} v_{\l}.
$$
The quantum binomials $\qbin{\l(h_i)+s}{n}{q_i}$ are
in the ring $\A$, and so it follows that
$\Lnum{K_i}{s}{n}\cdot v_{\l} \in \A v_{\l}$.  Similarly,
$\Lnum{D}{s}{n} \cdot v_{\l} \in \A v_{\l}$.  Hence,
$Z \cdot v_{\l} \in \A v_{\l}$.

The remainder of the proof now follows that of Theorem 3.5.  We
observe that the coefficients of terms of lower degree obtained
when commuting basis elements are in $\A$.  We noted in Section 1
that these coefficients are in fact in $\C[q,q^{-1}]$. The result
then follows from Theorem 3.5. \hskip 1cm \qed
\enddemo

Now that we have a vector space basis for the $\A$-form
$M_J^{\A}(\l)$ of $M_J^q(\l)$, we can begin comparing the two modules,
first as vector spaces.

\proclaim{Proposition 4.3} For any $\l \in P$ and $J \subseteq \dot I$,
as $\C(q)$-vector spaces,
$\C(q) \ot_{\A} M_J^{\A}(\l) \isom M_J^q(\l)$.
\endproclaim
\demo{Proof}The proof is fairly standard.
The $\C(q)$-linear map
$\zeta : \C(q) \ot_{\A} M_J^{\A}(\l) \to M_J^q(\l)$
defined by $\zeta(f \ot v) = fv$ for $f \in \C(q)$ and
$v \in M_J^{\A}(\l)$ is clearly surjective.  Let
$\{E_{\omega}\cdot v_{\l}\ |\ \omega\in \Omega\}$ be the basis
of $M_J^q(\l)$ determined by Theorem 3.5.  Let
$\xi :M_J^q(\l) \to \C(q) \ot_{\A} M_J^{\A}(\l)$ be a $\C(q)$-linear
map defined by
$$\xi(E_{\omega}\cdot v_{\l}) =
1 \ot E_{\omega}\cdot v_{\l}.
$$
Then, by  Proposition 4.2, $\xi$ is well-defined and the maps
$\zeta$ and $\xi$ are inverses.
\hskip 1cm \qed
\enddemo

We define a weight structure on $M_J^{\A}(\l)$ by setting
$M_J^{\A}(\l)_{\mu}= M_J^{\A}(\l) \cap M_J^q(\l)_{\mu}$
for each $\mu \in P$.

\proclaim{Proposition 4.4}
$M_J^{\A}(\l)$ is a weight module with the weight decomposition
$M_J^{\A}(\l) = \oplus_{\mu \in P}M_J^{\A}(\l)_{\mu}$.
\endproclaim
\demo{Proof}  The proof is quite standard, as in
\cite{BKMe, Proposition 3.23}.
\hskip 1 cm \qed
\enddemo

The vector-space isomorphism given above restricts to each weight
space and we obtain the following result.
\proclaim{Proposition 4.5}
For each $\mu \in P$, $M_J^{\A}(\l)_{\mu}$ is a free $\A$-module
and $\text{rank}_{\A}M_J^{\A}(\l)_{\mu}= \dim_{\C(q)}
M_J^q(\l)_{\mu}$.
\endproclaim

\b
\subhead 5. Classical limits \endsubhead
\b

In this section we take the classical limits of the $\A$-forms of
the quantum Verma-type modules, and show that
they are isomorphic to the Verma-type modules of $U(\g)$.

Recall that
$\A= \C[q,q^{-1}, \frac{1}{[n]_{q_i}}, i\in I, n>0]$.
Let
$\J$ be the ideal of $\A$ generated by $q-1$.  Then there
is an isomorphism of fields $\A / \J \isom \C$
given by $f + \J \mapsto f(1)$ for any $f \in \A$.
For any untwisted affine Kac-Moody algebra
$\g$, let $U_{\A} = U_{\A}(\g)$, and set
$U' = (\A/\J) \ot_{\A} U_{\A}$.
Then $U' \isom U_{\A}/\J U_{\A}$.  Denote by $u'$
the image in $U'$ of an element $u \in U_{\A}$.  It was shown
by Lusztig \cite{Lu} and DeConcini and Kac \cite{DK} that
$(D')^2=1$ and $(K_i')^2 =1$ for all $i \in I$.  If we let
$K'$ denote the ideal of $U'$ generated by $D'-1$ and
$\{K_i'-1\ |\ i \in I\}$, then $\o{U} = U'/K' \isom U(\g)$,
the universal enveloping algebra of $\g$.

Note that, under the natural map
$U_{\A} \to U_{\A}/\J U_{\A} \isom U'$, we have $q \mapsto 1$.
The composition of natural maps
$$
U_{\A} \to U_{\A}/\J U_{\A} \isom U' \to \o{U} = U'/K' \isom U(\g),
$$
is called taking the {\it classical limit} of $U_{\A}$.

Let $\o u \in \o U$ denote the image of an element $u \in U_{\A}$.
Then $\o U$ is generated by the elements $\o{E_i}, \o{F_i},
\o{D} := \o {\Lnum{D}{0}{1}}$ and
$\o{H_i}:= \o{ \Lnum{K_i}{0}{1}}$, $i \in I$, and, under
the isomorphism between $U(\g)$ and $\o U$, the elements
$e_i, f_i, d$ and $h_i$ may be identified with
$\o{E_i}, \o{F_i}, \o{D}$ and $\o {H_i}$, respectively.
Further, Beck \cite{Be1, Section 6} showed that
we may identify the $\o {E_{\beta}}$ with a PBW basis of
$U(\g)$, with elements denoted $e_{\beta}$.

\m
For $\l \in P$, $J \subseteq \dot I$, let
$M_J'(\l) = \A/\J \ot_{\A} M_J^{\A}(\l)$.  Then
$M_J'(\l) \isom M_J^{\A}(\l) / \J M_J^{\A}(\l)$ and $M_J'(\l)$ is
a $U'$-module.  For $\mu \in P$, let
$M_J'(\l)_{\mu} = \A / \J \otimes_{\A} M_J^{\A}(\l)_{\mu}$.
Since $M_J^{\A}(\l) = \oplus_{\mu \in P} M_J^{\A}(\l)_{\mu}$,
we must have $M_J'(\l) = \oplus_{\mu \in P} M_J'(\l)_{\mu}$.
 We also have the following standard result.

\proclaim{Proposition 5.1}
For $\mu \in P$, $\dim_{\A/\J} M_J'(\l)_{\mu} =
\text{rank}_{\A}M_J^{\A}(\l)_{\mu}$.
\endproclaim
\demo{Proof}
By Proposition 4.5, each weight space $M_J^{\A}(\l)_{\mu}$,
$\mu \in P$,
is a free $\A$-module.  Let $\{ v_j \ |\ j \in \Omega\}$ be a
basis for $M_J^{\A}(\l)_{\mu}$ .  Then every element
$v' \in M_J'(\l)_{\mu}= \A /\J \ot_{\A} M_J^{\A}(\l)_{\mu}$
can be written uniquely as
$v' = \sum_{j \in \Omega} a_j \ot v_j$ for some scalars
$a_j \in \A /\J$. (see \cite{Hu, Chapter 4, Theorem 5.11}).
Hence, $\{1 \ot v_j \ |\ j \in \Omega \}$ is a basis
for $M_J'(\l)_{\mu}$.
\hskip 1cm \qed
\enddemo

\proclaim{Proposition 5.2} The elements $D'$ and
$K_i'$ ($i \in I$) in $U'$ act as the identity on the $U'$ module
$M_J'(\l) = \A/\J \ot_{\A} M_J^{\A}(\l)$.
\endproclaim
\demo{Proof}
Let $\mu \in P$ and $\{ v_j\ |\ j \in \Omega\}$ be a basis of
$M_J^{\A}(\l)_{\mu}$.  Then by Proposition 5.1,
$\{ v_j' = 1\ot v_j \ |\ j \in \Omega\}$ is an $\A/\J$-basis
for $M_J'(\l)_{\mu}$.  Let $i \in I$.  For each $j \in \Omega$, we
have $K_i \cdot v_j = q_i^{\mu(h_i)}v_j$.  Letting $q \mapsto 1$,
we see $K_i' \cdot v_j' = v_j'$.  Thus, $K_i'$ acts on the identity
on each weight space of $M_J'(\l)$ and, since $M_J'(\l)$ is a weight
module, each $K_i'$ acts as the identity on the whole space.
Similarly, $D \cdot v_j = q_0^{\mu(d)}v_j$ implies that
$D' \cdot v'_j = v_j$, and that $D'$ acts as a scalar on $M_J'(\l)$.
\hskip 1cm \qed
\enddemo

Since $M_J'(\l)$ is a $U'$-module, $\o M_J(\l) = M_J'(\l)/K'M_J'(\l)$ is a
$\o U = U'/K'$-module.  But $K'$ was the ideal generated by
$D'-1$ and the $K_i'-1$, and $D'$ and each $K_i'$ acts as the
identity on $M_J'(\l)$, so $\o M_J(\l) = M_J'(\l)$.  Since $\o U \isom U(\g)$,
this means $\o M_J(\l)$ has a $U(\g)$-structure.  The module $\o M_J(\l)$
is called the classical limit of $M_J^{\A}(\l)$.  For $v \in M_J^{\A}(\l)$,
let $\o v$ denote the image of $v$ in $\o M_J(\l)$.

\proclaim{Proposition 5.3}
Let $v_{\l}$ be the generating vector for $M_J^{\A}(\l)$.  Then as a
$U(\g)$-module, $\o M_J(\l)$ is a weight module generated by
$\o{v_{\l}}$ and such that, for any $\mu \in P$, $\o M_J(\l)_{\mu}$ is
the $\mu$-weight space of $\o M_J(\l)$.
\endproclaim
\demo{Proof}
Let $v_{\l}$ generate $M_J^{\A}(\l)$, so that
$M_J^{\A}(\l) = U_{\A} \cdot v_{\l}$.  Then
$\o M_J(\l) = \o U \cdot \o {v_{\l}}$, so $\o{v_{\l}}$ generates
$\o M_J(\l)$.  As noted above, $M_J'(\l)$ is a $U'$-weight module
and since $\o M_J(\l) = M_J'(\l)$, $\o M_J(\l)$ is also a weight module.
Hence, $\o M_J(\l) = \oplus_{\mu \in P} \o M_J(\l)_{\mu}$.  It
remains to show that the vector space $\o M_J(\l)_{\mu}$ is
actually the $\mu$-weight space of $\o M_J(\l)$.  That is, we
have to show that $h_i \cdot \o{v_{\mu}} = \mu(h_i) \o {v_{\mu}}$
and $d \cdot \o{v_{\mu}} = \mu(d) \o {v_{\mu}}$ for all $i \in I$
and $\o {v_{\mu}} \in \o M_J(\l)_{\mu}$.

For $v_{\mu} \in M_J^{\A}(\l)_{\mu}$ and $i \in I$, we have
$$
\Lnum{K_i}{0}{1}\cdot v_{\mu} = \frac{K_i - K_i^{-1}}{q_i -
q_i^{-1}}\cdot v_{\mu}
= \frac{ q_i^{\mu(h_i)} - q_i^{-\mu(h_i)}}{q_i - q_i^{-1}} v_{\mu}.
$$
Passing to the classical limit, we obtain
$$
h_i \cdot \o{v_{\mu}} = \o{H_i} \cdot \o{v_{\mu}}
= \mu(h_i) \o{v_{\mu}}.
$$
Similarly, we have $d \cdot \o{v_{\mu}} = \o {D} \cdot \o {v_{\mu}}
= \mu(d) \o{v_{\mu}}$.
\hskip 1cm \qed
\enddemo

\proclaim{Proposition 5.4}
The $U(\g)$-module $\o M_J(\l)$ is a $J$-highest weight,
$U(\frak{n}_{-J})$-free module.
\endproclaim
\demo{Proof}
Let $v_{\l}$ be a generating vector for $M_J^q(\l)$.  By definition,
$E_{\beta} \cdot v_{\l} = 0$ for all $\beta \in S_J$.  Hence,
$E_{\beta} \cdot v_{\l} = 0$ in the $\A$-form $M_J^{\A}(\l)$.
Thus we have $\o{E_{\beta}}\cdot \o{v_{\l}}=0$ in $\o M_J(\l)$.
By Proposition 5.3, $\o M_J(\l)$ is a weight module generated
by $\o{v_{\l}}$ and $U(\frak{n}_{J})$ is spanned by the
$\o{E_{\beta}}$, $\beta \in S_J$, so $\o M_J(\l)$ is a $J$-highest
weight module.

It remains to show that $\o M_J(\l)$ is a free $U(\frak{n}_{-J})$-module.
>From Proposition 4.2, we know that $M_J^{\A}(\l)$ is isomorphic
to the space spanned by the ordered monomials
$E_{-\beta-n\d}\dots E_{-\beta+k\d}\dots
E_{-\a-n\d}\dots E_{-k\d}\dots E_{\a-k\d}$,  for
$\a \in \dot \D_+\cap \D^J$, $\beta \in \dot \D_+ \cap \D_J^{\infty}$,
$n\ge 0$, $k>0$.  Hence, $\o M_J(\l)$ is
isomorphic to the space spanned by the ordered monomials
$\o{E_{-\beta-n\d}}\dots \o{E_{-\beta+k\d}}\dots
\o{E_{-\a-n\d}}\dots \o{E_{-k\d}}\dots \o{E_{\a-k\d}}$.
But monomials in the images
$\o{E_{-\beta-n\d}}$,
$\o{E_{-\beta+k\d}}$, $\o{E_{-\a-n\d}}$, $\o{E_{-k\d}}$
and $\o{E_{\a - k\d}}$ form a basis
for $U(\frak{n}_{-J})$ and so $\o M_J(\l)$ is a
$U(\frak{n}_{-J})$-free module.
\hskip 1cm
\qed
\enddemo

We have shown that, for any $\l \in P$, if we start with a quantum
Verma-type module $M_J^q(\l)$, construct the $\A$-form
$M_J^{\A}(\l)$ and take the classical limit, then the resulting module
$\o M_J(\l)$ is a $U(\g)$-module isomorphic to the Verma-type
module $M_J(\l)$.  We have also seen that the weight space structure
 is preserved under these operations, and so $M_J^q(\l)$ is a true
quantum deformation of $M_J(\l)$.

\proclaim{Proposition 5.5} Let $\g$ be an affine Kac-Moody algebra.
Let $\l \in P$, $J \subseteq \dot I$.  Then the Verma-type
module $M_J(\l)$ admits a quantum deformation to the quantum
Verma-type module $M_J^q(\l)$ over $U_q(\g)$ in such a
way that the weight space decomposition is preserved.  In particular,
we have
$\dim M_J^q(\l)_{\mu} \neq 0$ if and only if $\l-\mu$ is in
the monoid generated by $S_J$,
$\dim M_J^q(\l)_{\l} = 1$, and
$0< \dim M_J^q(\l)_{\mu} < \infty$ if and only if $\l -\mu \in Q^J_+$.
\endproclaim

\m
\subhead 6. Submodule structure of quantum Verma-type modules
\endsubhead
\m

Using the quantum deformation results obtained above, we are
now in a position to prove some structural results about quantum
modules of Verma-type.

Let $\g_q^J$ be the subalgebra of $U_q(\g)$ generated by
monomials of the
form $X_i^{fin}$ and $Y_i^{fin}$ ($i = 1,2,3$) and $U_q^0(\g)$.
Let $\tilde \g_q^f$ denote the algebra generated by monomials of the
form $X_1^{fin}, X_3^{fin},Y_1^{fin}$ and $Y_3^{fin}$.

We conjecture that $\tilde \g_q^{f} \isom U_q(\g^{f})$, where
$U_q(\g^{f})$ is the subalgebra of $U_q(\g)$ generated by
$E_i, F_i$, ($i \in J$), $K_i$ ($i\in I$) and $D$, but we do not
have a proof in general.

Let
$G_q$ be the subalgebra of $U_q(\g)$ generated by
monomials of the form $X_2$ and $Y_2$.

 Then
the following proposition holds.

\proclaim{Proposition 6.1} Let $Z$ be the center of $G_q$.  Then
there exists a subalgebra $\o G_q$ of $G_q$ such that
$$ \g_q^J \isom \tilde \g_q^f\ot_Z \o G_q. $$
\endproclaim
\demo{Proof}
Recall the subalgebra $\o G$ of Proposition 2.2.  An element
$x \in \o G$ consists of $\C$-linear combinations of imaginary
root vectors $e_{k\d}^{(i)}$.  Write $x = \sum c_{i,k} e_{k\d}^{(i)}$.
For each $k \in \Z$ and $i=1, \dots, N$, by Beck's construction
there is a root vector $E_{k\d}^{(i)}$ in $U_q(\g)$ such that
the classical limit of $E_{k\d}^{(i)}$ may be identified with
$e_{k\d}^{(i)}$.  Consider the element $X = \sum c_{i,k} E_{k\d}^{(i)}$
in $U_q(\g)$.  Since the coefficients of the $E_{k\d}^{(i)}$'s
are in $\C$, certainly $X \in U_{\A}$ and the classical limit
of $X$ may be identified with $x \in \o G$.  Write $X = \phi(x)$.
Then we define a subalgebra $\o G_q$ of $G_q$ by
$$
\o G_q = \langle f(q)X \ |\ f(q) \in \C(q), X = \phi(x)
\text{ for some } x \in \o G\rangle.
$$

Let $\tilde G_q = G_q \cap \tilde \g_q^f$.  We now show that $G_q
\subseteq \tilde G_q \times \o G_q$.  Suppose $E_{k\d}^{(i)} \in
G_q$.  Then, by Proposition 2.2, we can write $e_{k\d}^{(i)} = e +
e'$, where $e \in \tilde \g^f$ and $e' \in \o G$.  Let $E' =
\phi(e')$ be the corresponding element in $\o G_q$.  Since $e \in
\tilde \g^f$, $e$ can be expressed as a linear combination of Lie
products of basis vectors $e_{\a+k\d}$ from root spaces $\pm \a +
n \d$, $\a \in \dot \D^J$. In $U_q(\g)$ the corresponding $E_{\a +
n \d}$'s are all in $\tilde \g_q^f$ and so $\C$-linear
combinations of products are in $\tilde \g_q^f$ as $\tilde \g_q^f$
is a subalgebra of $U_q(\g)$. Hence, there is some $E$ in $\tilde
\g_q^f$ such that $E_{k\d}^{(i)} = E + E'$.  This shows that all
the generators of $\g_q^J$ can be realized in the product $\tilde
\g_q^f \times \o G_q$.

Next we show that $[\tilde \g_q^f, \o G_q] =0$.
It is enough to show the $E_{\a_j +l\d}^{(i)}$ ($j \in J$)
commute with all $X = \phi(x)$, $x \in \o G$.  Assume
$X = \sum c_{i,k} E_{k\d}^{(i)}$.  Then
$$
\aligned
[E_{\a_j+l\d}, X] &= [E_{\a_j+l\d}, \sum c_{i,k} E_{k\d}^{(i)}] \\
&= \sum c_{i,k} [E_{\a_j+l\d}, E_{k\d}^{(i)}] \\
&= (\pm) \sum c_{i,k} \frac{1}{k}[k a_{ij}]_{d_i} E_{\a_j +
(k+l)\d},
\endaligned
$$
where we obtain the last equality by \cite{BK 1.6.5b}, and
neglect the sign.  All the coefficients of the last right-hand
side are in $\A$. Taking classical limits, we obtain
$$
[e_{\a_j+l\d}, x] = (\pm)\sum c_{i,k} a_{ij} e_{\a_j + (k+l)\d}.
$$
But $[e_{\a_j+l\d}, x]=0$, by Proposition 2.2, and $\sum c_{i,k}
a_{ij} e_{\a_j + (k+l)\d}=0$ if and only if $\sum c_{i,k}
\frac{1}{k}[k a_{ij}]_{d_i} E_{\a_j + (k+l)\d}=0$. Hence,
$[E_{\a_j+l\d}, X]=0$, and so $[\tilde \g_q^f, \o G_q] =0$.
Since $\tilde \g_q^f \cap \o G_q = Z$, we conclude that
$\g_q^J \isom \tilde \g_q^f \ot_Z \o G_q$.
\hskip 1cm \qed
\enddemo

Let $M^q_f(\l)$ be an ordinary Verma module for $\g^J_q$.
Since $\g^J_q\isom \tilde \g^f_q \ot_Z \o G_q$ by Proposition
6.1, we have $M^q_{f}(\l) \isom \tilde M^q_f(\l) \ot \o M^q(\l)$,
where $\tilde M^q_f(\l)$ is a ``Verma'' module for
$\tilde \g^f_q$ and
$\o M^q(\l)$ is a Verma module for $\o G_q$.  Further,
$\o M^q(\l)$ is
irreducible if and only if $\l(c)\neq 0$ as we show below.

Let $v_{\l}$ be a generating vector for the quantum Verma-type
module $M^q_J(\l)$ with $J$-highest weight $\l$.

\proclaim{Proposition 6.2} The $G_q$-module
$H^q(\l) = G_q\cdot v_{\l}$ is irreducible if and only if $\l(c) \neq 0$.
\endproclaim
\demo{Proof} Let $T_q=\{ E_{-\d}^{(i)} \cdot v_{\l}\}$, and set
$V_q = G_q \cdot T_q$.  Then $V_q \neq (0)$ and $V_q$ is a
$G_q$-submodule of $H^q(\l) = G_q \cdot v_{\l}$. Consider the
classical limits.  We have $\o T = \{ \o{E_{-\d}^{(i)}}\cdot \o
v_{\l}\} = \{ e_{-\d}^{(i)}\cdot \o v_{\l}\}$ and $\o V= \o G_q
\cdot \o T = G \cdot \o T$.  But $\o V$ is a non-zero submodule of
$H(\l)$.  If $\l(c) = 0$, $\o V$ is a proper submodule of $H(\l)$
and so $V_q$ is a proper submodule of $H^q(\l)$.  In
particular, $H^q(\l)$ is not irreducible.

On the other hand, if $\l(c) \neq 0$, then $H(\l)$ is
irreducible.  Suppose $W$ is a proper submodule of $H^q(\l)$.
Since $H^q(\l) = G_q\cdot v_{\l}$, using the construction of
Proposition 4.3, there is some module $W_{\A}$ such that
$W \isom \C(q)\ot_{\A}W_{\A}$.  We refer to this module as
the $\A$-form of $W$. Taking classical limits, we see that
$\o W_{\A}$ is  a proper submodule of $H(\l)$.  Hence,
we must have $H^q(\l)$ is irreducible if and only if $\l(c) \neq 0$.
\hskip 1cm \qed
\enddemo

\proclaim{Corollary 6.3} The $\o G_q$-module $\o M^q(\l)$
is irreducible if and only if $\l(c) \neq 0$.
\endproclaim

Let $U_q^{-J}$ be the $\C(q)$-linear subspace of $U_q(\g)$ spanned
by monomials of the form $Y_1^{\infty}$ and $X_3^{\infty}$.
It
follows from Theorem 3.5 that
$M_J^q(\l) \isom U_q^{-J} \ot M_f^q(\l)$ as vector spaces.

\proclaim{Theorem 6.4} Let $N$ be a non-trivial submodule of
$M_J^q(\l)$. Then
\newline
(i)  $N^{f} := N \cap M^q_f(\l) \neq 0$, and \newline (ii) if
$\l(c) \neq 0$, then as vector spaces, $N \isom U_q^{-J}
\ot N^f$.
\endproclaim
\demo{Proof} (i) We recall the basis of $M^q_J(\l)$ constructed in
Theorem 3.5.  Any element of $M^q_J(\l)$ may be written in the
form $x\cdot v_{\l} \in M^q_J(\l)$ where $x$ represents the sum of
ordered monomials in this basis. It is enough to consider
homogeneous $x$. Set $R(x\cdot v_{\l})$ to be minus the sum of the
heights of the finite roots not in $\dot \D^J$ in the
decomposition of $x\cdot v_{\l}$ (i.e., each $-\alpha+k\delta$
contributes $\text{ht}(\alpha)$) if $x\cdot v_{\l} \neq 0$ and set
$R(0) = 0$. Then it is clear that $R(x\cdot v_{\l})=0$ if and only
if $x \in \g^J_q$.

It is enough to show that there exists $y\in U_q(\g)$ such that
$yx\cdot v_{\l}\neq 0$ and $R(yx\cdot v_{\l})< R(x\cdot v_{\l})$.
We will find an element $y=E_{\a-K\d}$ where $K$ is sufficiently
large and $\a$ is some suitable root.

Let $w=x \cdot v_{\l} \in M_J^q(\l)$.  Then, as in Proposition 4.3, we
may write $w=fw'$ for some $f \in \C(q)$ and $w' \in M_J^{\A}(\l)$.
Then $w' = f^{-1}w=f^{-1}x\cdot v_{\l} \in U_q\cdot w$.
Furthermore, suppose $w' = (q-1)^kw''$, with $k>0$ and
$w'' \in M_J^{\A}(\l)$.  Then $w'' = (q-1)^{-k}w' \in U_q\cdot w$.
Hence, without loss of generality, we may assume
$w=x\cdot v_{\l} \in M_J^{\A}(\l)$ and $q-1$ is not a factor of $w$ in
$M_J^{\A}(\l)$.  Taking the classical limit, we then have
$\o w = \o{x \cdot v_{\l}} \neq0$.

Suppose $\o w$ is in $M_J^f(\l)$.  Since $x$ is homogeneous, the
grading of $M_J^q(\l)$ ensures that $w$ is in $M_f^q(\l)$. Suppose
that $\o w$ is not in $M_J^f(\l)$. Then by \cite{Fu5, Lemma 5.4}
there exists a root $\a$ and nonnegative integer $K$ such that
$e_{\a-K\d}\o w\neq 0$. Note that $e_{\a-K\d}$ is the image of
$E_{\a-K\d}$.  Hence, $E_{\a-K\d}x\cdot v_{\l} \neq 0$ and
$R(E_{\a-K\d}x\cdot v_{\l} )<R(x\cdot v_{\l})$. We complete the
proof by induction.  This proves part (i) of the theorem.

Let $0 \neq w \in N$.  As in the proof of part (i), we can assume
that $w \in M_J^{\A}(\l)$.  Using the basis of $M_J^q(\l)$, we may
write $w = \sum u_i' u_i \cdot v_{\l}$, where $u_i' \in
U_q^{-J}$ and $u_i \in \frak{n}^q_{- J}$ for each $i$. It
follows that $\o u_i'\o u_i \cdot \o v_{\l} \neq 0$ for each $i$.

We will assume that each $u'_i$ is a monomial of type
 $$ u_i' =
E_{-\beta_{i,1}+n_{i,1}\d}^{l_{i,1}}\dots
E_{-\beta_{s(i),i}+n_{s(i),i}\d}^{l_{s(i),i}}, $$ where
$\beta_{i,j} \in \dot \D_+ \setminus \D^J$, and
$-\beta_{i,j}+n_{i,j}\d \neq -\beta_{k,j}+n_{k,j}\d$ for $i \neq
k$ and $u_i' \neq u_j'$ for $i \neq j$.  Further, without loss of
generality, we may also assume $w$ is homogeneous.

Then $\o w = \sum \o u_i' \o u_i \cdot \o v_{\l}$ and $$ \o u_i' =
e_{-\beta_{i,1}+n_{i,1}\d}^{l_{i,1}}\dots
e_{-\beta_{s(i),i}+n_{s(i),i}\d}^{l_{s(i),i}}.$$ By
\cite{Fu5, Theorem 5.14(i)}, we conclude that $\o u_i \cdot v_{\l}
\in \o N^f$ for all $i$. Hence, $u_i \cdot v_{\l}$ belongs to
$N_f^{\A}$, the $\A$-form of $N$, and so $u_i \cdot v_{\l}$ is in
$N_f$ for each $i$.  The statement follows. \hskip 1cm \qed
\enddemo

\proclaim{Corollary 6.5} $M_J^q(\l)$ is irreducible if and only if
$M_f^q(\l)$ is irreducible as $\g_q^J$-module if and only if
$\l(c) \neq 0$ and $\tilde M^q_f(\l)$ is irreducible as
$\tilde \g_q^f$-module.
\endproclaim
\demo{Proof}  The results follows from Corollary 6.3, Theorem 6.4(i),
and the fact that $M_f^q(\l) \isom \tilde M^q_f(\l) \ot
 \o M^q(\l)$.
\enddemo

\m
\subhead 7. Imaginary Verma modules
\endsubhead
\m

In the special case of imaginary Verma modules we can say rather
more about their structure.  In particular, we can give a precise
description of their submodule structure in the case of level
zero.  Throughout this section, we shall suppose that
$J=\emptyset$ so that the modules $M^q_{\emptyset}(\l)$ and
$M_{\emptyset}(\l)$ are imaginary Verma modules for $U_q(\g)$ and
$U(\g)$, respectively.  We shall also
assume that $\l(c)=0$. In this case $H^q(\l)$ is a
reducible $G_q$-module with maximal submodule $H^q_0(\l)$
consisting of all spaces except $\C v_{\l}$.

Denote $M^q_{\emptyset,0}(\l)=U_q (\g)H^q_0(\l)$.  We remark that
$M^q_{\emptyset,0}(\l) \neq M_\emptyset^q(\l)$.

Set $\widetilde{M_\emptyset^q(\l)} = M_\emptyset^q(\l)/M^q_{\emptyset, 0}(\l)$.

\proclaim{Theorem 7.1}The $U_q(\g)$-module
$\widetilde{M_\emptyset^q(\l)}$ is irreducible if and only if
 $\l(h_i)\neq 0$, for all $i=1,\dots, N$.
\endproclaim
\demo{Proof} Let $\widetilde{M_\emptyset^q(\l)}$ be irreducible and suppose
there exists some $i \in \{1, \dots, N\}$ such that $\l(h_i)=0$.
Set $\a = \a_i$ and $E = E_{-\a}$. We have that
$\widetilde{M_\emptyset^q(\l)}=U_q(\g)\cdot v_{\l}$.
Consider $W=U_qE\cdot v_{\l}$. We show $W$ is a proper
submodule of  $\widetilde{M_\emptyset^q(\l)}$.

Since $E \cdot v_{\l}\neq 0$,
$W \neq (0)$.  Suppose $W = \widetilde{M_\emptyset^q(\l)}$.
Then there must exist elements $m_j$ in $U_q$
such that $v_{\l}=\sum_j m_j E\cdot v_{\l}$, and for each $m_j$,
we have $\text{ht}(m_j)=\a$.  (Recall the height of a
monomial is the sum of heights of finite roots involved.)
Using Beck's ordering and the notation introduced in the proof
of Theorem 3.5, we may write each $m_j$ in the form
$m_j = Y_1Y_2Y_3ZX_3X_2X_1$.  We consider the actions of the
$m_j$ on $E\cdot v_{\l}$.  We need the following lemma.

\proclaim{Lemma 7.2}For any $k \in \Z$ and $\beta \in \dot \Delta_+$,
$E_{-\beta+k\d}E\cdot v_{\l}=0$.
\endproclaim
\demo{Proof}We divide the proof into cases.\newline
1. $k\in \Z$, $\beta \neq \a$.

 Using Beck's basis, we write
$E_{\beta+k\d}E \cdot v_{\l}= \sum_i n_i^-n_i^+v_{\l}$
where $n_i^+$ are ordered monomials of the form $X_3X_2X_1$ and
$n_i^-$ are ordered monomials of the form $Y_1Y_2Y_3$.
Since $\beta \neq \a$ and $\a$ is simple, then, due to the convexity
of the basis, this ordered expression must
contain on the right an element of the form $X_1$ or $Y_3$
(depending on the sign of $k$).  But
$X_1\cdot v_{\l}= Y_3 \cdot v_{\l} = 0$. Hence
 $E_{\beta+k\d}E \cdot v_{\l}= 0$.
\newline
2. $k \in \Z \setminus\{ 0\}$, $\beta = \a$.

By \cite{BK, (1.6.5d)},
$E_{\a + k\delta}E\cdot v_{\l} = X_{k\delta}\cdot v_{\l}=0$
for some suitable vector $X_{k\delta}$ of weight $k\delta$,
and this has a trivial action in $\widetilde{M_\emptyset^q(\l)}$.
\newline
3. $k=0$, $\beta = \a$.

In this case $E_{\a}E_{-\a} =0$ because $\l(h_i)=0$.
\newline
In all cases, we have $E_{-\beta+k\d}E\cdot v_{\l}=0$.
\hskip 1 cm \qed
\enddemo

\ni{\bf Return to Proof of Theorem}

Consider the action of an element $Y_1Y_2Y_3ZX_3X_2X_1E\cdot v_{\l}$.
By Lemma 7.2, $X_1E\cdot v_{\l}=0$.  Now $X_2 \cdot v_{\l} = 0$, and
so by \cite{BK, (1.6.5b)}, $X_2E \cdot v_{\l}$ is in the space of
elements of the form $X_3 \cdot v_{\l}$.  Elements of the form $Z$
commute with elements of the form $X_3$ and act as scalars on $v_{\l}$.
As shown in Theorem 3.5, we may commute elements of the form $Y_3$
with elements of the form $X_3$ and then $Y_3 \cdot v_{\l} = 0$.
Hence we are left considering elements of the form
$Y_1Y_2X_3 \cdot v_{\l}$.  The monomial is non-zero as it contains
$E$ and has height 0 as we supposed
$\sum m_j E \cdot v_{\l} = v_{\l}$.  This is not
possible and we have a contradiction.  Hence $W$ is a proper submodule
of $\widetilde {M_\emptyset^q(\l)}$ and $\widetilde{M_\emptyset^q(\l)}$ is reducible.

Now we prove the converse statement of the theorem.
Let $\l(h_i)\neq 0$, $i=1,\dots, N$.
Let $\widetilde{M_\emptyset^{\A}(\l)} =
\widetilde{M_\emptyset^q(\l)}\cap M_\emptyset^{\A}(\l)$
denote the $\A$-form of $\widetilde{M_\emptyset^q(\l)}$, and let
$\widetilde{M_\emptyset(\l)}$ denote its classical limit.  Suppose $N^q$ is
a proper submodule of $\widetilde{M_\emptyset^q(\l)}$.  Then
$N^{\A} = N^q \cap \widetilde{M_\emptyset^{\A}(\l)}$ is a submodule
of $\widetilde{M_\emptyset^{\A}(\l)}$ and the classical limit of $N^{\A}$
gives a proper submodule of $\widetilde{M_\emptyset(\l)}$.  But this last
module is irreducible by \cite{Fu3}.  The theorem is proved.
\hskip 1cm \qed
\enddemo

\proclaim{Proposition 7.3}
Let $\l(c)=0$ and $\l(h_i)\neq 0$ for all $i$. Then $M_\emptyset^q(\l)$
has an infinite filtration with irreducible quotients of the form
$\widetilde{M_\emptyset^q(\l+k\d)}$,
$k \geq 0$.
\endproclaim
\demo{Proof} Follows from theorem above.
\hskip 1cm \qed
\enddemo

\proclaim{Corollary 7.4}
 Let $\l(c)=0$, $\l(h_i)\neq 0$, $i=1,\dots,N$ and $N^q$ be a submodule
of $M_\emptyset^q(\l)$. Then $N^q$ is generated by $N^q\cap H^q(\l)$.
\endproclaim
\demo{Proof} Obvious from Proposition 7.3.
\hskip 1cm \qed
\enddemo

\b


\Refs
\widestnumber\key{XXXX}
\m

\ref \key Be1
\by J. Beck
\paper Braid group action and quantum affine algebras
\jour Commun. Math. Phys.
\vol 165
\yr 1994
\pages 555--568
\endref

\ref\key Be2
\bysame
\paper Convex bases of PBW type for quantum affine algebras
\jour Commun. Math. Phys.
\vol 165
\yr 1994
\pages 193--199
\endref

\ref\key BK
\by J. Beck and V.G. Kac
\paper Finite-dimensional representations of quantum affine
algebras at roots of unity
\jour J. Amer. Math. Soc.
\vol 9
\yr 1996
\pages 391--423
\endref

\ref\key BKMe
\by G. Benkart S.-J. Kang and D.J. Melville
\paper Quantized enveloping algebras for Borcherds superalgebras
\jour Trans. Amer. Math. Soc.
\vol 350
\yr 1998
\pages 3297--3319
\endref

\ref\key CP
\by V. Chari and A. Pressley
\book A Guide to Quantum groups
\publ Cambridge University Press
\publaddr Cambridge
\yr 1994
\endref

\ref\key Co1
\by B. Cox
\paper Verma modules induced from nonstandard Borel subalgebras
\jour Pac. J. Math.
\vol 165
\yr 1994
\pages 269--294
\endref

\ref\key Co2
\bysame
\paper Structure of nonstandard category of highest weight modules
\inbook Modern Trends in Lie algebra representation theory
\bookinfo V. Futorny, D. Pollack (eds.)
\pages 35--47
\yr 1994
\endref


\ref\key CFKM
\by B. Cox, V.M. Futorny, S.-J. Kang and D.J. Melville
\paper Quantum deformations of imaginary Verma modules
\jour Proc. London Math. Soc.
\vol 74
\yr 1997
\pages 52--80
\endref

\ref\key CFM
\by B. Cox, V.M. Futorny and D.J. Melville
\paper Categories of nonstandard highest weight modules for affine Lie algebras
\jour Math. Z.
\vol 221
\yr 1996
\pages 193--209
\endref

\ref\key Da
\by I. Damiani
\paper La $R$-matrice pour les alg\`ebres quantiques de type
affine non tordu
\jour preprint
\endref

\ref \key DK
\by C. DeConcini and V.G. Kac
\paper Representations of quantum groups at roots of 1
\inbook Operator Algebras, Unitary Representations, Enveloping
Algebras and Invariant Theory
\bookinfo A.Connes, M.Duflo, A.Joseph, R.Rentschler (eds)
\publ Birkh\"auser
\publaddr Boston
\pages 471--506
\yr 1990
\endref

\ref\key Dr
\by V.G. Drinfel'd
\paper Hopf algebras and the quantum Yang-Baxter equation
\jour Soviet Math. Dokl.
\vol 32
\yr 1985
\pages 254--258
\endref

\ref\key Fu1
\by V.M. Futorny
\paper Root systems, representations and geometries
\jour Ac. Sci. Ukrain. Math.
\vol 8
\yr 1990
\pages 30--39
\endref

\ref\key Fu2
\bysame
\paper The parabolic subsets of root systems and corresponding
representations of affine Lie algebras
\jour Contemp. Math.
\vol 131
\yr 1992
\pages 45--52
\endref

\ref \key Fu3
\bysame
\paper Imaginary Verma modules for affine Lie algebras
\jour Canad. Math. Bull.
\vol 37
\yr 1994
\pages 213--218
\endref

\ref\key Fu4
\bysame
\paper Verma type modules of level zero for affine Lie algebras
\jour Trans. Amer. Math. Soc.
\vol 349
\yr 1997
\pages 2663--2685
\endref

\ref\key Fu5
\bysame
\book Representations of affine Lie algebras
\bookinfo Queen's Papers in Pure and Applied Mathematics, vol 106
\publ Queen's University
\publaddr Kingston
\yr 1997
\endref

\ref\key FGM
\by V.M. Futorny, A.N. Grishkov and D.J. Melville
\paper Quantum imaginary Verma modules for affine Lie algebras
\jour C.R. Math. Rep. Acad. Sci. Canada
\vol 20
\yr 1998
\pages 119--123
\endref

\ref\key Ga
\by F. Gavarini
\paper A PBW basis for Lusztig's form on untwisted affine quantum groups
\jour preprint
\vol q-alg$\backslash$9709018
\yr 1997
\endref

\ref\key Hu
\by T. Hungerford
\book Algebra, (5th edition)
\publ Springer-Verlag
\publaddr New York
\yr 1989
\endref

\ref\key JK1
\by H.P. Jakobsen and V.G. Kac
\paper A new class of unitarizable highest weight representations
of infinite dimensional Lie algebras
\inbook Lecture Notes in Physics
\vol 226
\publ Springer
\publaddr Berlin
\yr 1985
\pages 1--20
\endref

\ref\key JK2
\bysame
\paper A new class of unitarizable highest weight representations of infinite di
mensional Lie algebras II
\jour J. Funct. Anal.
\vol 82
\yr 1989
\pages 69--90
\endref

\ref \key Ja
\by J.C. Jantzen
\book Lectures on Quantum Groups
\publ American Mathematical Society
\publaddr Providence
\yr 1996
\endref

\ref\key Ji
\by M. Jimbo
\paper A $q$-difference analogue of $U(\g)$ and the Yang-Baxter equation
\jour Lett. Math. Phys.
\vol 10
\yr 1985
\pages 63--69
\endref


\ref \key K
\by V.G. Kac
\book Infinite-dimensional Lie algebras (3rd edition)
\publ Cambridge University Press
\publaddr Cambridge
\yr 1990
\endref

\ref\key Ka
\by S.-J. Kang
\paper Quantum deformations of generalized Kac-Moody algebras and
their modules
\jour J. Algebra
\vol 175
\yr 1995
\pages 1041--1066
\endref

\ref\key KT
\by S.M. Khoroshkin and V.N. Tolstoy
\paper On Drinfeld's realization of quantum affine algebras
\jour J. Geom. Phys.
\vol 11
\yr 1993
\pages 445--452
\endref

\ref\key Le
\by F.W. Lemire
\paper Note on weight spaces of irreducible linear representations
\jour Canad. Math. Bull.
\vol 11
\yr 1968
\pages 399--403
\endref

\ref \key Lu
\by G. Lusztig
\paper Quantum deformations of certain simple modules over enveloping algebras
\jour Adv. Math.
\vol 70
\yr 1988
\pages 237--249
\endref

\ref\key M
\by D.J. Melville
\paper An $\A$-form technique of quantum deformations
\inbook Recent developments in quantum affine algebras and related
topics, Contemp. Math., vol 248
\publ Amer. Math. Soc.
\yr 1999
\pages 359--375
\endref

\endRefs
\enddocument